\documentclass{article}[12pt]
\usepackage[utf8]{inputenc}
\usepackage[margin=1in]{geometry}
\usepackage{xcolor}
\usepackage{amsmath,amssymb,amsthm} 
\usepackage{hyperref}
\usepackage{pgfplots}
\usepackage{graphicx,epstopdf,caption}

\usepackage{braket}
\usepackage{amsfonts}
\usepackage{color}
\usepackage{amsmath}
\usepackage{booktabs}
\usepackage{siunitx}
\usepackage{caption}

\usepackage{xcolor}
\usepackage{tikz}
\usetikzlibrary{shapes,arrows}
\tikzstyle{block} = [rectangle, draw, fill=blue!20, text width=10em, text centered, rounded corners, minimum height=4em]
\tikzstyle{line} = [draw, -latex']

\usepackage{pgfplots}
\pgfplotsset{filter discard warning=false}
\pgfplotsset{compat=1.18}
\pgfplotsset{
  plot style/.style={every axis plot/.append style={#1}},
  bar style/.style={every axis plot post/.append style={#1}},
  axis style/.style={every axis/.append style={#1}},
  mark style/.style={every mark/.append style={#1}},
  numformat/.style={/pgf/number format/#1},
  xnumformat/.style={xticklabel style={numformat/.list={#1}}},
  ynumformat/.style={yticklabel style={numformat/.list={#1}}},
  lognumformat/.style={log plot exponent style/.append style=
    {numformat/.list={#1}}},
  scatter classes/.style={only marks,scatter,
    scatter src=explicit symbolic,scatter/classes={#1}},
  axis style={
    xlabel near ticks, ylabel near ticks, % Better label placement
    every axis title shift=3pt, % Better title placement
    scale only axis, % Natural axis sizes
    cycle list={}, % Disable cycle list
  },
  plot style={mark style={solid}},
  table/col sep=comma, % Use comma separated value data tables
  table/search path={./data},
}

\usepackage{fixme}
\fxsetup{
  status=draft,
  theme=color,
  silent,
}
\fxsetup{marginface=\linespread{1}\footnotesize}
\FXRegisterAuthor{tk}{tke}{\colorbox{orange!50}{\color{orange}TK}}
\FXRegisterAuthor{jk}{jke}{\colorbox{olive!50}{\color{olive}JK}}
\FXRegisterAuthor{rw}{rwe}{\colorbox{purple!50}{\color{purple}JK}}
\FXRegisterAuthor{rj}{rje}{\colorbox{teal!50}{\color{teal}JK}}

\usetikzlibrary{matrix}

\usepackage[capitalise,nameinlink,noabbrev]{cleveref}
\newtheorem{theorem}{Theorem}[section]
\newtheorem{example}[theorem]{Example}
\newtheorem{lemma}[theorem]{Lemma}

\newtheorem{proposition}[theorem]{Proposition}

\newtheorem{remark}{Remark}
\newtheorem{definition}{Definition}

\crefname{lemma}{Lemma}{lemmas}
\crefname{proposition}{Proposition}{propositions}
\crefname{definition}{Definition}{definitions}
\usepackage{algorithm}
\usepackage{algpseudocode}
\newenvironment{varalgorithm}[1]
  {\algorithm\renewcommand{\thealgorithm}{#1}}
  {\endalgorithm}

\usepackage{subcaption}
\usepackage{array}
\usepackage{pgfplots}
\usepackage{graphicx,epstopdf}
\usepackage{tikz}
\usepackage[mathscr]{euscript}
\usepackage{nicefrac}
\usepackage{xspace}

\def\mL{\mathbf{L}}
\def\mH{\mathbf{H}}
\def\mG{\mathbf{G}}

\def\mP{\mathbf{P}}
\def\mR{\mathbf{R}}

\def\mQ{\mathbf{Q}}
\def\Real{\mathbb{R}}
\def\mH{\mathbf{H}}
\def\mM{\mathbf{M}}
\def\mA{\mathbf{A}}
\def\mB{\mathbf{B}}

\def\mC{\mathbf{C}}
\def\mS{\mathbf{S}}

\def\nor{\mathcal{N}}
\def\vx{\mathbf{x}}

\def\v{\mathbf{v}}
\def\vy{\mathbf{y}}
\def\vz{\mathbf{z}}
\def\vw{\mathbf{w}}
\def\vf{\mathbf{f}}

\def\E{\mathbb{E}}
\def\EU{\text{EU}}

\def\eps{\pmb{\varepsilon}}

\def\mX{\mathbf{X}}
\def\mY{\mathbf{Y}}
\def\mZ{\mathbf{Z}}
\def\mM{\mathbf{M}}
\def\mO{\mathbf{O}}
\def\S{\mathbb{S}}
\def\qr{\texttt{orth}}
\def\O{\mathcal{O}}

\def\tMs{\pmb{\mathscr{M}}^d}

\def\tX{\pmb{\mathscr{X}}}
\def\tY{\pmb{\mathscr{Y}}}

\def\tr{\text{trace}}

\def\id{\mathbf{I}}
\def\tC{\pmb{\mathscr{C}}}
\def\tM{\pmb{\mathscr{M}}}
\def\mL{\mathbf{L}}
\def\vw{\mathbf{w}}

\def\gr{\operatorname{Gr}}
\def\st{\operatorname{St}}

\def\diag{\text{diag}}
\def\mat{\operatorname{mat}}
\def\SNR{\texttt{SNR}}

\def\step{\pmb{\gamma}}
\def\mR{\mathbf{R}}
\def\grad{\operatorname{grad}}
\def\hess{\operatorname{Hess}}
\def\tangent{\text{T}}
\def\colspan{\text{colspan}}
\def\ad{\text{ad}}
\def\D{\text{D}}
\def\sym{\operatorname{sym}}

\def\qrcode{\text{qr}}

\title{Scalable symmetric Tucker tensor decomposition}
\author{Ruhui Jin\thanks{Department of Mathematics, University of Wisconsin-Madison, Madison, WI} \and
Joe Kileel\thanks{Department of Mathematics, University of Texas at Austin, Austin, TX} \and
Tamara G. Kolda\thanks{MathSci.ai, Dublin, CA}
\and and Rachel Ward\footnotemark[2] }
\date{\vspace{-3ex}}

\begin{document}

\maketitle
\begin{abstract}
We study the best low-rank Tucker decomposition of symmetric tensors.
The motivating application is decomposing higher-order multivariate moments. Moment tensors have special structure and are important to various data science problems.
We advocate for projected gradient descent (PGD) method and higher-order eigenvalue decomposition (HOEVD) approximation as computation schemes. 
Most importantly, we develop scalable adaptations of the basic PGD and HOEVD methods to decompose sample moment tensors.
With the help of implicit and streaming techniques, we evade the overhead cost of building and storing the moment tensor.
Such reductions make computing the Tucker decomposition realizable for large data instances in high dimensions. 
Numerical experiments demonstrate the efficiency of the algorithms and the  applicability of moment tensor decompositions to real-world datasets.
Finally we study the convergence on the Grassmannian manifold, and prove that the update sequence derived by the PGD solver achieves first-
and second-order criticality. 
\end{abstract}

\section{Introduction}
\label{sec: introduction}
\begin{figure}[!htb]
  \centering
  \usetikzlibrary{3d}
  \begin{tikzpicture}[thick, join=bevel, scale=0.8,
      mat/.style n args={2}{transform shape, draw, anchor=north west,minimum height=#1cm, minimum width=#2cm}
    ]
    % Draw tensor X
    \draw[join=bevel]
    (0,0) coordinate (X1SW) -- ++(0,2,0) -- ++(0,0,-2) -- ++(2,0,0) coordinate (X2NE) -- ++(0,-2,0) -- ++(0,0,2) -- cycle
    (0,2,0) -- ++(2,0,0) coordinate (X1NE) -- ++(0,0,-2)
    (2,2,0) -- (2,0,0);
    % Draw equal sign
    \path (X2NE) ++(0,-1) node[right,font=\large,inner sep=0.5em] (eq) {$\approx$};
    % Draw Q matrix on the left
    \node[mat={2}{1}] at (eq.east |- X1NE) (Q1) {};
    \path (Q1.north east) ++(0.8em,0) ++(0,-1) coordinate (CC);
    % Draw tensor C
    \begin{scope}[shift=(CC)]
      \draw[join=bevel]
      (0,0) coordinate (C1SW) -- ++(0,1,0) -- ++(0,0,-1) coordinate (C2NW) -- ++(1,0,0) coordinate (C2NE) -- ++(0,-1,0) -- ++(0,0,1) -- cycle
      (0,1,0) -- ++(1,0,0) coordinate (C1NE) -- ++(0,0,-1)
      (1,1,0) -- (1,0,0);
    \end{scope}
    % Draw Q matrix on the right
    \path (C2NE |- C1SW) ++(0.8em,0) node[mat={1}{2},anchor=south west] (Q2) {};
    % Draw Q matrix on top
    \begin{scope}[shift=(C2NW)]
      \path[canvas is xz plane at y=0] ++(0,-1.5em) node[mat={2}{1}] (Q3) {};
    \end{scope}
    % Add labels
    \begin{scope}
      \path (X1SW) -- (X1NE) node[midway] {$\tX$};
      \path (C1SW) -- (C1NE) node[midway] {$\tC$};
      \node at (Q1.center) {$\mQ$};    
      \node at (Q2.center) {$\mQ$};    
      \node at (Q3.center) {$\mQ$};          
    \end{scope}
  \end{tikzpicture}
  \caption{Symmetric Tucker decomposition.} 
  \label{fig: sym Tucker}
\end{figure}
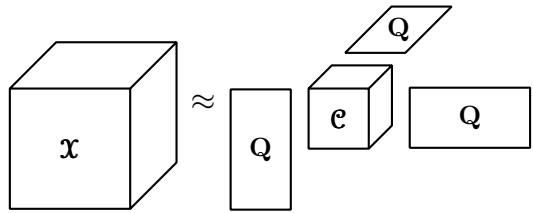

Tensor decompositions \cite{KB09} provide powerful tools to represent and compress higher-order arrays.
In this work, we focus on the \emph{symmetric rank-$r$ Tucker decomposition}. % \cite{T66}.  
Applied to a $d$-th order symmetric tensor of size $n \times \ldots \times n$, it expresses the input using a 
symmetric core tensor $\tC\in\Real^{r \times \dots \times r}$  ($r \ll n$) and a  matrix $\mQ \in \Real^{n \times r}$ with orthonormal columns. 
See \cref{fig: sym Tucker} for  an \nolinebreak illustration.

%%% Local Variables:
%%% mode: latex
%%% TeX-master: "fig1_standalone"
%%% End:

\emph{Symmetric tensors} abound in statistics as a means of representing higher-order moments and cumulants of a random vector. 
They are commonly used to analyze dependences in multivariate distributions \cite{C94,HO00,LMV002}, learn latent variable models via the method-of-moments \cite{ARDSM14,GHK15}, solve important signal processing problems \cite{C91,C98}, etc.

In most cases, the exact moment of the ground-truth distribution is not available.
The idea of moment estimation is to use the $d$-th order \emph{sample moment tensor}, denoted $\tMs$, instead.  
This is formed from observations $\{\vx_1, \ldots, \vx_p\} \subseteq \Real^n$ of a random variable $\vx$ as 
\begin{equation}
\label{eqn: sample moment}
\tMs = \frac{1}{p}\sum_{i=1}^p \vx_i^{\otimes d} \approx \E(\vx^{\otimes d}).
\end{equation}
There is a need for economical representations of $\tMs$, which would otherwise take $\O(n^d)$ space to store. 
Using the symmetric Tucker decomposition, one approximates $\tMs$ solely by determining the rank-$r$ factor matrix $\mQ$, which takes just $rn$ memory space. (The core $\tC$ is a byproduct once $\mQ$ is identified; see \cref{rem:core}.)
However, decomposing the sample moment still has upfront costs of $\O(pn^d)$ and $\O(n^d)$ to initially construct and store $\tMs$.
For large data volumes $p$ and high dimensionalities $n$, we turn to \emph{implicit} (or \emph{tensor-free}) and \emph{streaming} strategies to escape such prohibitive overhead.

The symmetric Tucker decomposition of moment tensors has other important motivations, besides compression.
One major application is to extract the key features $\mQ \in \Real^{n \times r}$ from moment tensors. 
Without the need for the core tensor $\tC \in \S^d(\Real^r),$ this usage further reduces the output's storage.
In some scientific domains, data has high volatility, e.g., in combustion simulations \cite{A19}, outlier detection \cite{GSJZ14,RHRS22} and asset allocation \cite{ML09,BKJ07,BCV20}.  Higher-order moments and cumulants ($d \geq 3$) are sensitive to heavy tails and anomalies. 
Thus, finding the key subspace $\mQ$ of variability in higher-order statistics is important for non-Gaussian (also known as non-normal) data analysis.
This is akin to the role of principal component analysis (PCA) and the covariance matrix for Gaussian data. 

In other settings, data observations are assumed to lie in a low-dimensional subspace and follow the \textit{low-rank factor model} \cite{S1987}. For example, this holds in the Fama-French model \cite{FF93} of finance and also in psychometrics \cite{M85}. 
The resulting moment tensors naturally admit a Tucker format. 
Thus the decomposition enforces the low Tucker rank structure and improves moment estimation accuracy.

\subsection{Related works}
\label{subsec: related works}
Due to the NP-hard nature of tensor decompositions \cite{HL13}, an approximate symmetric Tucker decomposition is computed by higher-order eigenvector decomposition (HOEVD) \cite{LMV001}. 
It is the eigendecomposition of the flattened tensor. 
It can be viewed as the symmetric version of higher-order singular value decomposition (HOSVD) \cite{LMV002} and efficient variants \cite{VVM12,MB18,SGLTU20,MSK20,CWY20,randhosvd21}. 
Despite their popularity, HOEVD offers no guarantee of returning the exact critical points for the problem.

We instead turn to iterative search methods to achieve the best possible solution, using HOEVD as an initialization. 
Newton \cite{ES09,ILAH09}, quasi-Newton \cite{SL10} and trust-region \cite{IAHL11} methods have been considered for general Tucker decomposition. These second-order optimizations rely on Riemannian Hessian updates on a the manifold of orthonormal matrices.  For symmetric tensors, higher-order orthogonal iterations (HOOI) \cite{LMV003}, Jacobi rotations \cite{IAV13, LUC18} and power method \cite{R13} were proposed. 
Regalia briefly mentions a shifted variant of the power method in \cite{R13}, which is akin to a gradient descent method, although strong convergence analysis is omitted in the work.

Our work is primarily motivated by moment tensor decompositions, which are fundamental objects in statistics and other disciplines. Neither the existing HOEVD nor iterative algorithms for Symmetric Tucker take advantage of the sample moment tensor structure \eqref{eqn: sample moment}. 
We show that streaming and implicit techniques can be applied to address the computational burden in the data-rich and high-dimensional regime. 
For covariance matrices which are centered second-order moments, this is relevant to the streaming PCA scheme, also known as Oja's method \cite{O82}, and its adaptations \cite{MCJ13,YHW18,HW19}.

Recently, a similar strategy was applied to optimization with moment tensor CP decomposition for learning Gaussian mixture models, see  \cite{SK20,PKK22}.  
Building on techniques from \cite{SK20}, the present work also passes over only a subset of samples at a time and drives the optimization without forming gigantic moment tensors.
Despite using similar approaches for increasing efficiency, we emphasize that CP and Tucker decompositions are useful for their own independent purposes. The motivations we list for feature extraction and data factor analysis specifically apply to Tucker decompositions.  We refer readers to \cref{subsec: Tucker vs CP} for a concrete comparison.

\subsection{Our contributions}
 In this paper, our major contribution is developing scalable algorithms for symmetric Tucker decomposition on moment tensors.
 
To start with, we propose the projected gradient descent (PGD) framework to the general symmetric Tucker tensor decomposition problem, see \cref{basic PGD}. 
The approach is led by a first-order gradient update governed by the Tucker approximation cost. It is free of non-trivial second-order updates on manifolds. 

Driven by their importance in data science, our attention then focuses on decomposing sample moment tensors \eqref{eqn: sample moment}.
 By exploiting the special input structure, we develop scalable PGD (SPGD) in \cref{algorithm}.  There  are two main highlights.
First, we adopt the implicit technique and operate on the input sample vectors directly. This much cheaper in memory and computation than constructing the sample moment tensor. Proposed in \cite{SK20,PKK22}, the implicit strategy greatly reduces the effort for gradient computation.
Second, SPGD is built on a streaming model. 
By passing over small batches of data observations,
both storage and computational complexity are further reduced. Similar ideas have been used in CP moment tensor decomposition \cite{SK20} for separate purposes.

We also design a scalable implementation of HOEVD (SHOEVD) for sample moment tensors
\eqref{eqn: sample moment}, see \cref{alg: scalable hoevd}. 
This method is akin to the streaming PCA idea. 
We use the SHOEVD approximation (\cref{alg: scalable hoevd}) as the initialization to SPGD (\cref{algorithm}). Moreover, SHOEVD may be of independent interest due to its broad usages, e.g. \cite{L97,AG05}.
The savings from the streaming and implicit strategies allow SHOEVD and SPGD to handle large-scale computations. The gains are summarized in \cref{cost}. Here $p$ is the total sample number and $b$ is the batch size at each iteration. Usually, the user sets $b \ll p$. These efficiency advantages are corroborated by extensive numerical experiments.  

We demonstrate the key motivations of sample moment decompositions, including feature extraction and moment tensor low-rank estimation. Our proposed SPGD is respectively applied to real datasets in anomaly detection and portfolio allocation.

In terms of theoretical guarantee, the basic PGD (\cref{basic PGD}) convergence result is provided. Employing the theories from manifold optimization and dynamical systems, we prove that for any symmetric tensor, PGD converges monotonically to a first-order critical point given arbitrary initialization.  Furthermore, it converges to a second-order critical point almost surely.  
 
\begin{table}[!htb]
{\footnotesize
  \caption{Computation and storage of SPGD (\cref{algorithm}) and SHOEVD (\cref{alg: scalable hoevd}).
    Here $p$ is the number of samples, $n$ is the vector size, $d$ is the moment tensor order, and $b$ is the batch size.}
\label{cost}
\begin{center}
  \begin{tabular}{ccc} 
  \toprule
 & \bf overhead  & \bf storage per iter.  
 \\[1mm]
\bf streaming $\&$ implicit & none & $nr+nb$  
\\[1mm]
\bf full $\&$ implicit & none & $nr+np$ 
\\ [1mm]
\bf full $\&$ explicit & $n^dp$ &$nr+n^d$  
\\ 
 \bottomrule[1pt]
  \end{tabular}
\end{center}
}
\end{table}

%%% Local Variables:
%%% mode: latex
%%% TeX-master: "../manuscript"
%%% End:

\section{Symmetric Tucker tensor decomposition}
\label{sec: sym Tucker}
\subsection{Notation and background}
Scalars, vectors, matrices and tensors are represented as lowercase letters $x$, lowercase bold letters $\vx$, uppercase bold letters $\mX,$ and calligraphic bold letters $\tX$, respectively.
The $(i_1, \dots, i_d)$-th element of a $d$-way tensor is denoted by $x_{i_1, \dots, i_d}$. The $j$-th column of a matrix $\mX$ is $\vx_j$. 
The spectral norm of a matrix is  $\|\cdot\|_2$. The norm of a tensor is $\| \cdot\|$, see \cref{def: inner product}.
The identity matrix is $\id_n \in \Real^{n \times n}$. 
The Moore-Penrose pseudoinverse of a matrix is denoted using $\dag$. 
The matrix symmetrization operator is written $\sym: \Real^{n \times n} \to \S^2(\Real^n)$ and given by $\sym(\mX) = \tfrac{1}{2} (\mX+ \mX^\top)$ for $\mX \in \Real^{n \times n}$.
We write $\qr$ for the function that selects the orthonormal factor $\mQ \in \Real^{n \times r}$ from the QR decomposition of a full-rank $n \times r$ matrix.\footnote{For ease of analysis, we ensure the uniqueness of the QR factorization by using the one where all diagonal entries of the $\mR$ factor are positive.  But none of our algorithms require this choice.} 
The notation $\mX^{[d]}$ is the elementwise exponentiation, where each coordinate of $\mX$ is raised to the $d$-th power.
The set $\{1, \ldots, n\}$ is denoted by $[n]$ when $n$ is a positive integer. 

\begin{definition}[symmetric tensors]
A tensor $\tX \in \Real^{n \times \dots \times n}$ is symmetric if its entries are invariant under any permutation of indices, i.e.
\[
x_{i_1, \dots, i_d} = x_{i_{\sigma(1)}, \dots, i_{\sigma(d)}},~~\forall (i_1, \dots, i_d) \in [n]^d~\text{and}~\sigma \in \Pi(d)
\]
where $\Pi(d)$ is the permutation group on $[d]$. 
Here, $d$ and $n$ are called the order and dimension of $\tX$.
The space of all real-valued symmetric tensors of order $d$ and dimension $n$ is denoted by $\S^d(\Real^n)$. 
\end{definition}

\begin{definition}[symmetric outer product]
The outer product of a vector $\vx \in \Real^n$ with itself $d$ times is denoted as $\vx^{\otimes d} \in \S^d(\Real^n)$. The $(i_1, \dots, i_d)$-th entry of this product is
\[
\left(\vx^{\otimes d}\right)_{i_1, \dots, i_d} = x_{i_1}\!\cdots x_{i_d},~~\forall (i_1, \dots, i_d) \in [n]^d. 
\]
\end{definition}

\begin{definition}[symmetric Tucker product]
\label{sym Tucker}
For a symmetric tensor $\tX \in \S^d(\Real^n)$ and a matrix $\mY\in\Real^{n \times m}$, their symmetric Tucker product is defined as $\tX$ with $\mY$ multiplied along each mode, written
$
\tX \cdot \left(\mY, \dots, \mY \right) \in\S^d(\Real^m).
$
The $(j_1, \dots, j_d)$-th entry of the product is
\[
\left(\tX \cdot \left(\mY, \dots, \mY \right)\right)_{j_1, \dots, j_d} = \sum_{i_1=1}^n\cdots \sum_{i_d=1}^n x_{i_1, \dots, i_d} \, y_{i_1, j_1} \dots y_{i_d, j_d}, ~~\forall (j_1, \dots, j_d) \in [m]^d.
\]
\end{definition}

{\sloppy 
\begin{definition}[symmetric tensor matricization]
A symmetric tensor $\tX \in \S^d(\Real^n)$ can be flattened into a matrix, denoted $\mat(\tX) \in \Real^{n \times n^{d-1}}$. 
The index $(i_1, \dots, i_d) \in [n]^d$ of $\tX$ is mapped to the index $\left(i_1,1+ \sum_{\ell = 2}^d (i_\ell -1)\,n^{\ell-2}\right)\in [n] \times [n^{d-1}]$ of $\mat(\tX)$.
\end{definition}}

\begin{definition}[inner product]
\label{def: inner product}
For (not necessarily symmetric) tensors $\tX, \tY \in \Real^{n_1\times n_2 \dots \times n_d}$ of the same size, their inner product is 
\[
\langle \tX, \tY\rangle = \sum_{i_d =1}^{n_d}\dots \sum_{i_1=1}^{n_1} x_{i_1,\dots,i_d} y_{i_1,\dots,i_d} \in \Real.
\]
The norm of a tensor is $\left\| \tX\right\|= \sqrt{\langle \tX, \tX \rangle}.$
\end{definition}

\begin{definition}
For tensors $\tX \in \Real^{m_1\times n_2 \dots \times n_d}$ and $\tY \in \Real^{m'_1\times n_2 \dots \times n_d},$ their inner product in all modes except the first one \cite{ES09,SL08} is defined as: 
\begin{multline*}
\mZ = \langle \tX, \tY\rangle_{-1} \in \Real^{m_1 \times m'_1}\\
\text{where}~~z_{j_1, j'_1} = \sum_{i_d=1}^{n_d}\dots\sum_{i_2=1}^{n_2} x_{j_1,i_2, \dots, i_d}y_{j'_1,i_2, \dots, i_d}~~~ \forall j_1 \in [m_1], j'_1 \in [m'_1].
\end{multline*}
\end{definition}

\subsection{Problem statement} \label{subsec:prob-state}

Let $\tX\in \S^d(\Real^n)$ be a symmetric tensor.  Fix a user-specified  rank $r$ (typically $r \ll n$).  Our goal is to find the best approximation of $\tX$ represented by the symmetric Tucker product of a core tensor $\tC \in \S^d(\Real^r)$ and a rank-$r$ orthonormal matrix $\mQ\in\Real^{n \times r}$. The problem is formulated as follows:
\begin{equation}
\label{eqn: sym Tucker C}
\min_{\tC\in\S^d(\Real^r) \atop \mQ \in \Real^{n \times r}} \left\| \tX - \tC \cdot \left(\mQ^\top, \dots, \mQ^\top \right) \right\|^2,~~\text{subject~to}~\mQ^\top\mQ=\id_r.
\end{equation}
In fact, the minimization \eqref{eqn: sym Tucker C} is equivalent to solving for a rank-$r$ basis $\mQ$ that maximizes the following cost function:
\begin{equation}
\tag{$*$}
\label{eqn: sym Tucker}
\max_{\mQ^\top\mQ = \id_r} F(\mQ) \equiv \left\| \tX \cdot \left(\mQ, \dots, \mQ\right)\right\|^2. 
\end{equation}

\begin{remark} \label{rem:core}
The problem formulation transforms from \eqref{eqn: sym Tucker C} to \eqref{eqn: sym Tucker} by solving the unconstrained linear least squares in \eqref{eqn: sym Tucker C} for the core, i.e., 
\begin{equation}
\label{eqn: core}
\tC=\tX \cdot \left(\mQ, \dots, \mQ \right).
\end{equation}
Substituting \eqref{eqn: core} into \eqref{eqn: sym Tucker C} gives
\begin{equation*}
\begin{array}{rl}
  \left\| \tX - \tC \cdot \left(\mQ^\top, \dots, \mQ^\top \right) \right\|^2
  & = \left\| \tX- \tX \cdot \left(\mQ\mQ^\top, \dots, \mQ\mQ^\top \right) \right\|^2\\
  & = \left\| \tX\right\|^2 -  \left\| \tX \cdot \left(\mQ\mQ^\top, \dots, \mQ\mQ^\top \right) \right\|^2\\
  & = \left\| \tX\right\|^2 -  \left\| \tX \cdot \left(\mQ, \dots, \mQ \right) \right\|^2.
\end{array}
\end{equation*}
As $\|\tX\|^2$ is a constant in $\mQ$, minimizing the above quantity is inded the same as maximizing $\left\| \tX \cdot \left(\mQ, \dots, \mQ \right) \right\|^2$.
\end{remark}

\begin{remark}
The full Tucker decomposition \eqref{eqn: sym Tucker C} takes $\mathcal{O}(nr + r^d)$ in memory to store the basis $\mQ\in\Real^{n \times r}$ and the core $\tC \in \S^d(\Real^r).$
Though still exponential in $d$, this cost savings can be substantial compared to the original $\tX \in \S^d(\Real^n)$ requiring $\mathcal{O}(n^d)$ space if the rank $r$ is much less than the ambient dimension $n$. 
Further, in applications such as tensor PCA \cite{MR14} and anomaly detection \cref{subsec: anomaly detection}, one only desires the feature subspace $\mQ$, and the core $\tC$ need note be computed or stored.  This reduces the cost to $\mathcal{O}(nr)$.
\end{remark}

The feasible domain in \eqref{eqn: sym Tucker} for $\mQ$ is the Stiefel manifold:
\begin{equation}
\label{eqn: stiefel}
\st(n,r): = \left\{\mQ \in \Real^{n\times r} ~\vert~ \mQ^\top\mQ = \id_r\right\},
\end{equation}
which is a compact Riemannian submanifold of $\Real^{n\times r}$ with dimension $\frac{r(2n-r-1)}{2}$. 
Further, the cost function $F$ is rotationally invariant. 
In fact, it only depends on the orthogonal projector $\mP = \mQ\mQ^\top \in \S^2(\Real^n)$, since $\left\| \tX \cdot \left(\mQ, \dots, \mQ \right)\right\| = \left\| \tX \cdot \left(\mP, \dots, \mP \right)\right\|$. 
So, the optimization \eqref{eqn: sym Tucker} naturally runs on the Grassmannian manifold, which consists of all orthogonal projectors $\mP $:
\begin{equation}
\label{eqn: grassmann}
\gr(n,r): = \left\{\mP \in \S^2(\Real^n) ~\vert~ \mP^2 = \mP, ~\text{rank}(\mP) = r\right\}.
\end{equation}
Here, $\gr(n,r)$ is a compact Riemannian manifold of $\S^2(\Real^n)$ with dimension $r(n-r)$. 

\begin{remark} \label{rem:eckart-young}
No simple characterization of the critical points of \eqref{eqn: sym Tucker} is known when $d \geq 3$.
This is in contrast to case of $d=2$, where the classical Eckart-Young theorem \cite{EY36} holds. 
Various possible extensions of Eckart-Young to tensors have been shown to fail, e.g., see \cite{K03}.  
\end{remark}

\subsection{PGD algorithm}
\label{subsec: PGD}
This paper develops a projected gradient descent (PGD) method to solve the constrained optimization problem \eqref{eqn: sym Tucker}. The PGD framework works as follows: given an initial point, we move the current iterate $\mQ_{t-1}$ along the negative gradient direction by a step size of $\gamma_t > 0$ where the gradient is: 
\begin{equation}
\label{eqn: euclid grad}
\nabla F(\mQ_{t-1}) =2d\left\langle \tX\cdot\left( \id_n, \mQ_{t-1}, \dots ,\mQ_{t-1}\right),  \tX\cdot\left( \mQ_{t-1}, \dots ,\mQ_{t-1}\right)\right\rangle_{-1} \in \Real^{n \times r}.
\end{equation}
Then we project the update back to the Stiefel manifold \eqref{eqn: stiefel} by selecting the $\mQ$ factor from the QR decomposition. The procedure is described in \cref{basic PGD}.  See \cref{sec: analysis} for a convergence analysis of PGD.
\begin{varalgorithm}{1}
%\begin{algorithm}
\caption{PGD for symmetric Tucker decomposition}
 \textbf{Input:} symmetric tensor $\tX \in \S^d(\Real^n)\\$
number of iteration T, step sizes $\{\gamma_t\} \subseteq \Real$\\
initial guess $\mQ_{0}\in \st(n,r)$ \Comment{initialization}\\
 \textbf{Output:}  rank-$r$ orthonormal basis $\mQ \in \st(n,r)$ 
 \vspace{1mm}
 \begin{algorithmic}[1]
\For{$t = 1, \dots, T$}
\State $\mQ_t \gets \mQ_{t-1}+2d\gamma_t\,\left\langle \tX\cdot\left( \id_n, \mQ_{t-1}, \dots ,\mQ_{t-1}\right),  \tX\cdot\left( \mQ_{t-1}, \dots ,\mQ_{t-1}\right)\right\rangle_{-1}$ \\
\Comment{gradient update}
\State $\mQ_t \gets \qrcode(\mQ_t,0)$ \Comment{retraction}
\EndFor\\
\Return $\mQ \gets {\mQ}_{T
}$
\end{algorithmic}
\label{basic PGD}
%\end{algorithm}
\end{varalgorithm}

\subsection{HOEVD initialization} 
\label{subsec: hoevd}
An approximate solution to the best symmetric Tucker decomposition \eqref{eqn: sym Tucker} is the HOEVD \cite{LMV001} result:
\begin{equation}
\label{eqn: hoevd}
\displaystyle \mQ_{\text{hoevd}} \in \st(n,r) ~\text{given~by~the~} \text{leading~}r~\text{eigenvectors~of~}\mat(\tX)\mat(\tX)^\top.
\end{equation}
HOEVD is widely used for Tucker approximation, beause it can be computed via a direct matrix eigendecomposition. 

Generally speaking, the HOEVD solution is not a critical point of \eqref{eqn: sym Tucker} \cite{LMV001}. 
Its approximation error to a $d$-th order tensor $\tX$ is bounded \cite[Theorems 6.9-10]{H14} as
\[
\left\|\tX - \hat{\tX}_{\text{hoevd}} \right\| \leq \sqrt{d}\left\|\tX - \hat{\tX}_{\text{best}} \right\|.
\]
Here, $\hat{\tX}_{\text{hoevd}} = \tX\cdot(\mQ_{\text{hoevd}}\mQ_{\text{hoevd}}^\top, \dots, \mQ_{\text{hoevd}}\mQ_{\text{hoevd}}^\top)$ as in \cref{rem:core} and $\hat{\tX}_{\text{best}}$ is the optimum of \eqref{eqn: sym Tucker C}. 
We observe the suboptimality of the HOEVD solution in \cref{subsec: rank and noise}. 
However according to \cite[section 3.4]{LMV002}, the HOEVD solution usually belongs to an ``attractive" region around the desired local optima in numerical experiments.  
So it is a good candidate for initialization.

%%% Local Variables:
%%% mode: latex
%%% TeX-master: "../manuscript"
%%% End:

\section{Scalable sample moment tensor decomposition}
\label{sec: scalable sample moment decomposition}

The main interest in this work is the  decomposition of the sample moment tensor $\tMs \in \S^d(\Real^n)$ \eqref{eqn: sample moment}. 
We design scalable algorithms for the symmetric Tucker decomposition \eqref{eqn: sym Tucker}, and the HOEVD approximation \eqref{eqn: hoevd}, of $\tMs$. 
An average of symmetric outer products, the sample moment structure \eqref{eqn: sample moment} is the key to the scalability of the decomposition computation. Tensor operations with outer products can be reduced to much cheaper computations involving only their factor vectors. 
Specifically we use the following identities:
for vectors $\vx,\vy \in \Real^n, \vz \in \Real^m$ and a matrix $\mY \in \Real^{n \times m},$ it holds \cite{H14}
\begin{equation}
\label{eqn: Tucker outer product}
\vx^{\otimes d} \cdot \left(\mY, \dots, \mY \right) = \left(\mY^\top \vx\right)^{\otimes d}\in\S^d(\Real^m),
\end{equation}
\begin{equation}
\label{eqn: implicit 1}
\langle \vx^{\otimes d}, \vy^{\otimes d}\rangle = \langle \vx,\vy \rangle^d \in \Real, \quad\text{and}
\end{equation}
\begin{equation}
\label{eqn: implicit 2}
\langle \vz\otimes \vx^{\otimes(d-1)}, \vy^{\otimes d}\rangle_{-1}
%= \mat \! \left(\vz\otimes \vx^{\otimes(d-1)}\right)\mat \! \left(\vy^{\otimes d}\right)^{\! \top}
\! = \vz\, \langle \vx,\vy \rangle^{d-1} \vy^\top \in \Real^{m \times n}.
\end{equation}

\subsection{Scalable PGD algorithm}
\label{subsec: scalable PGD}

In \cref{algorithm}, we introduce the scalable PGD (SPGD) when the input tensor is a sample moment \eqref{eqn: sample moment}.
The symmetric Tucker cost \eqref{eqn: sym Tucker} now depends solely on the data observations $\mX = [\vx_1, \dots, \vx_p] \in \Real^{n\times p}.$ 
Indeed due to \eqref{eqn: Tucker outer product},

\begin{equation}
\label{eqn: sym Tucker moment}
F(\mQ) := \left\|  \left(\frac{1}{p}\sum_{i=1}^p \vx_i^{\otimes d}\right) \cdot \left(\mQ, \dots, \mQ\right)\right\|^2= \left\| \frac{1}{p}\sum_{i=1}^p \left(\mQ^\top \vx_i\right)^{\otimes d}\right\|^2. 
\end{equation}
We actually do not need to create or store the sample moment $\tMs$.
Such an idea was called the implicit technique in \cite{SK20,PKK22}. 

Commonly, there is a huge amount of data in exceedingly high dimensions. We hence further turn to the streaming strategy to speed up the computations.  
Here data is revealed as a sequence of sample batches $\{\mX_1, \dots, \mX_t, \dots \} \subseteq \Real^{n \times b}$ with batch size $b$ ($b \ll p$). 
We read each batch as it arrives but discard it afterwards. This streaming model is called the turnstile model \cite{M05,LNW14,TYUC19}.

Therefore, as the $t$-th $(t \geq 1)$ batch $\mX_t = [\vx_{t,1}, \dots, \vx_{t,b}]^\top \in \Real^{n \times b}$ arrives, the cost function \eqref{eqn: sym Tucker moment} is recast stochastically:
\begin{equation}
\label{eqn: stream cost}
F_t(\mQ):= \left\| \frac{1}{b}\sum_{i=1}^b \left(\mQ^\top \vx_{t,i}\right)^{\otimes d}\right\|^2.
\end{equation}
In light of \eqref{eqn: implicit 2}, we can compute the stochastic gradient in a tensor-free manner: 
\begin{equation}
\label{eqn: implicit gradient}
\begin{array}{ll}
\nabla F_t(\mQ) & \displaystyle = 2d\,\langle\frac{1}{b} \sum_{i=1}^b \vx_{t,i}\otimes \left(\mQ^\top \vx_{t,i}\right)^{\otimes (d-1)} , \frac{1}{b}\sum_{j=1}^b \left(\mQ^\top \vx_{t,j}\right)^{\otimes d}\rangle_{-1}\\
&\displaystyle=\frac{2d}{b^2}\, \sum_{i=1}^b\sum_{j=1}^b \vx_{t,i} \langle \mQ^\top \vx_{t,i}, \mQ^\top \vx_{t,j}\rangle^{d-1} \vx_{t,j}^\top \mQ\\
&\displaystyle = \frac{2d}{b^2}\,\mX_t \left(\mX_t^\top \mQ\mQ^\top \mX_t\right)^{[d-1]} \mX_t^\top \mQ \,\, \in \,\, \Real^{n \times r}.
\end{array}
\end{equation}
Here, $\left(~ \cdot ~\right)^{[d-1]}$ means raising each element of the input matrix to the $(d-1)$-th power. 
Computing the gradient implicitly via online samples reduces the cost to $\O(rnb+rb^2)$ flops per iteration. This improves upon $\O(rnp+rp^2)$ using the full amount of samples, and a cost that is exponential in $d$ 
if we operate on the explicit tensor in $\S^d(\Real^n)$.

Moreover, a scalable HOEVD solver, $\textsc{shoevd}$ in line 1 of \cref{algorithm}, is utilized for good initialization; see \cref{subsec: scalable HOEVD} for details. 
We also use AdaGrad \cite{DHS11} to automatically tune step sizes. 
We implement the column-wise extension of AdaGrad proposed in \cite{HW19}, since it is natural to independently update the orthonormal basis vectors in $\mQ$. The step size $(\step)_t$ becomes a row vector in $\Real^{1 \times r},$ rather than a scalar.
See lines 4 and 5, where $.  \slash$ denotes columnwise division.

SPGD is a two-phased iterative algorithm.  It consists of 
\begin{itemize}
\item
Phase I - SHOEVD approximation (line 1),
\item
Phase II - symmetric Tucker decomposition (lines 2-8). 
\end{itemize}
%The alterations with respect to \cref{basic PGD} are highlighted in blue below.
\begin{varalgorithm}{2}
%\begin{algorithm}
\caption{SPGD for sample moment Tucker decomposition}
 \textbf{Input:} 
 numbers of iterations $T_1, T_2$; batch sizes $b_1, b_2$; \\
 datastreams $\{\mX_1,\dots, \mX_{T_1}\} \in \Real^{n \times b_1}$ (Phase I), $\{\mX_{T_1{+}1}, \dots, \mX_{T_2}\} \in \Real^{n \times b_2} $ (Phase II); \\
hyperparameters: $(\step_1)_0 = (\step_2)_0 = 10^{-5}\text{ones}(1,r)$; $c_1, c_2 >0$ \\
 \textbf{Output:} rank-$r$ orthonormal basis $\mQ \in \st(n,r)$
 \vspace{1mm}
 \begin{algorithmic}[1]
\State $\mQ_{T_1} \gets \textsc{shoevd}( T_1, b_1, \{\mX_1,\dots, \mX_{T_1}\}, (\step_1)_0, c_1)$ \Comment{{\color{black}\textup{{initialization}}}}
\For{$t = T_1+1, \dots, T_2$}
\State $\mG_t \gets \frac{2d}{b_2^2}\,\mX_t \left(\mX_t^\top \mQ_{t-1}\mQ_{t-1}^\top \mX_t\right)^{[d-1]} \mX_t^\top \mQ_{t-1} $ \Comment{{\color{black}\textup{gradient computation}}}
\State $(\step_2)_t \gets \left((\step_2)_{t-1}^{[2]} + \text{colnorm}(\mG_t)^{[2]}\right)^{[1/2]}$ \Comment{{\color{black}\textup{step size}}}
\State $\mQ_t \gets \mQ_{t-1}+c_2\,\mG_t \, .  \slash  (\step_2)_t$
\State $\mQ_t \gets \qrcode(\mQ_t,0)$ 
\EndFor
\State \Return $\mQ \gets {\mQ}_{T_2
}$
\end{algorithmic}
\label{algorithm}
%\end{algorithm}
\end{varalgorithm}

%\tkwarning[layout=inline]{I think the square root in line 4 should be written using the bracket exponent notation to make it clear that it's elementwise.}

\subsection{Scalable HOEVD}
\label{subsec: scalable HOEVD}
We present the scalable HOEVD (SHOEVD) method in \cref{alg: scalable hoevd}. 
This serves as the Phase I  initialization in SPGD (\cref{algorithm}), but it may be of independent interest.
As HOEVD is essentially a matrix eigendecomposition, we follow the streaming PCA idea.
The solver uses the PGD framework, but is applied to a different cost function than \eqref{eqn: sym Tucker moment}.  Specifically, we solve HOEVD iteratively finding $\mQ \in \st(n,r)$ that maximizes the following stochastic cost formed by a data batch $\mX_t \in \Real^{n \times b_1}$:
\[
\left\| \mat \left(\frac{1}{b_1}\sum_{i=1}^{b_1} \vx_{t,i}^{\otimes d}\right)^\top\!\mQ\right\|^2\!.
\]
Using \eqref{eqn: implicit 2}, the gradient can be calculated implicitly: 
\begin{equation}
\label{eqn: hoevd grad}
\begin{array}{ll}
\displaystyle2\,\mat \left(\frac{1}{b_1}\sum_{i=1}^{b_1} \vx_{t,i}^{\otimes d}\right)\mat \left(\frac{1}{b_1}\sum_{j=1}^{b_1} \vx_{t,j}^{\otimes d}\right)^\top\! \mQ 
&\displaystyle=\frac{2}{b_1^2}\, \sum_{i=1}^{b_1}\sum_{j=1}^{b_1} \vx_{t,i} \langle \vx_{t,i}, \vx_{t,j}\rangle^{d-1} \vx_{t,j}^\top \mQ\\
&\displaystyle= \frac{2}{b_1^2}\,\mX_t\left(\mX_t^\top \mX_t\right)^{[d-1]} \mX_t^\top \mQ \in \Real^{n \times r}.
\end{array}
\end{equation}
It costs $\O(rnb_1+nb_1^2+rb_1^2)$ flops to compute the gradient.  
%See the whole procedures in \cref{alg: scalable hoevd}.

\begin{varalgorithm}{3}
%\begin{algorithm}
\caption{SHOEVD for sample moment decomposition}
 \textbf{Input:}
 number of iterations $T_1$; batch size $b_1$;\\
 data stream $\{\mX_1,\dots, \mX_{T_1}\} \subseteq \Real^{n \times b_1}$  ;
\\
hyperparameters: $(\step_1)_0 = 10^{-5}\text{ones}(1,r)$; $c_1 >0$ \\
 \textbf{Output:}  top-$r$ eigenvectors $\mQ_{\text{hoevd}} \in \st(n,r)$
 \vspace{1mm}
 \begin{algorithmic}[1]
 \Function{s-hoevd}{$T_1, \{\mX_1,\dots, \mX_{T_1} \}, (\step_1)_0, c_1$} \\
 $\mQ_{0} \gets [\mathcal{N}(0,1)]^{n \times r}$ \\ $\mQ_{0} \gets \qrcode(\mQ_{0},0)$ \Comment{initialization}
\For{$t = 1, \dots, T_1$}
\State $\mG_t \gets \frac{2}{b_1^2}\,\mX_t\left(\mX_t^\top \mX_t\right)^{[d-1]} \mX_t^\top \mQ_{t-1}$ 
\State $(\step_1)_t \gets \left( (\step_1)_t^{[2]}+ \text{colnorm}(\mG_t)^{[2]}\right)^{[1/2]}$ 
\State $\mQ_t \gets \mQ_{t-1}+c_1\,\mG_t \, .\slash (\step_1)_t $ \Comment{gradient update}
\State $\mQ_t \gets \qrcode(\mQ_t,0)$ \Comment{retraction}
\EndFor\\
\Return $\mQ_{\text{hoevd}} \gets \mQ_{T_1} $
\EndFunction
\end{algorithmic}
\label{alg: scalable hoevd}
%\end{algorithm}
\end{varalgorithm}

%\tkwarning[layout=inline]{I think the square root in line 6 should be written using the bracket exponent notation to make it clear that it's elementwise.}

\begin{remark}
SPGD and SHOEVD (\cref{algorithm,alg: scalable hoevd}) obviate the upfront cost to form the tensor, which would require $\O(pn^d)$ flops. 
Each iteration temporarily stores an updated matrix $\mQ_t \in \Real^{n \times r}$ and an online sample batch $\mX_t \in \Real^{n \times b}$, taking only $nr+nb$ memory.
This explains the complexity counts shown in \cref{cost}. 
\end{remark}

\begin{remark} 
\label{rem:streaming core}
The core tensor $\tC\in \S^d(\Real^r)$ \eqref{eqn: core} can be evaluated as follows.
After obtaining $\mQ \in \st(n,r)$ by SPGD or SHOEVD (\cref{algorithm,alg: scalable hoevd}), we draw fresh data $\mX_{\text{c}} = [(\vx_{\text{c}})_{1}, \dots, (\vx_{\text{c}})_{p_{\text{c}}}]^\top \in \Real^{n \times p_{\text{c}}}$ and then \nolinebreak set:
\[
\tC =  \frac{1}{p_{\text{c}}}\sum_{i=1}^{p_{\text{c}}} \left(\mQ^\top (\vx_{\text{c}})_{i}\right)^{\otimes d}.
\]
The core computation costs $\mathcal{O}(r^d)$ in storage and $\mathcal{O}(p_crn + p_cr^d)$ in flops.
\end{remark}

%%% Local Variables:
%%% mode: latex
%%% TeX-master: "../manuscript"
%%% End:

\section{Synthetic numerical tests}
\label{sec: numerics}
\label{subsec: synthetic}
In this section, we utilize synthetic tests to study the proposed algorithms for sample moment tensor decomposition. And real-world datasets are considered in \cref{sec: applications}. 
In this section, we show:
\begin{itemize}
    \item the scalability of implicit compared to explicit methods in \cref{imp_exp_time},
    \item the speedup of streaming compared to full data usage in \cref{str_full_time},
    \item the robustness of the HOEVD approximation and the symmetric Tucker decomposition to noise and target rank choice in \cref{noise,rank}, and
    \item a comparison between symmetric Tucker and CP in \cref{table: err 4}.
\end{itemize}
\subsection{Set-up}
\label{subsec: experimental set-up}
In these tests, we generate moment tensors with approximate low-rank symmetric Tucker structure \eqref{sym Tucker}. One way is to assume data instances $\vx \in \Real^n$ follow a linear factor model:
\begin{equation}
\label{eqn: factor model}
\vx = \mB\vf +\pmb{\eps}.
\end{equation}
Here we fix the factor loading matrix $\mB \in \Real^{n \times r_{\text{true}}}$. 
The latent factor $\vf \in \Real^{r_{\text{true}}}$ is drawn from a coordinate-independent, zero-mean and unit variance normal-inverse Gaussian (NIG) distribution \cite{B97,BCV20}. This distribution is skewed and heavy-tailed, aligning with some motivations discussed in \cref{sec: introduction}. 
The noise ${\pmb{\eps}} \in \Real^n$ is Gaussian, following $\nor(\mathbf{0}, \sigma^2\id_n)$. 

Note that in \eqref{eqn: factor model}, the $d$-th moment of $\vx$ roughly admits a low-rank Tucker format:
\begin{equation*}
\begin{array}{lll}
\tMs \!\!&\approx\E\left(\vx^{\otimes d}\right)\!\!&\approx \,\,\, \displaystyle\E\left((\mB\vf)^{\otimes d}\right) \\
&& \displaystyle= \,\,\,  \E\left((\mR^*\vf)^{\otimes d}\right)\cdot \left(\mQ^{*\top}, \dots, \mQ^{*\top}\right).
\end{array}
\end{equation*}
The matrices $\mQ^*\in \Real^{n\times r_{\text{true}}}, \mR^* \in \Real^{r_{\text{true}} \times r_{\text{true}}}$ are from the QR factorization of $\mB$. 

We define the inverse signal-to-noise ratio
\begin{equation}
\label{eqn: SNR}
\SNR^{-1} := \frac{\| \E\left(\pmb{\eps}\right)\|}{\| \E\left(\mB\vf\right)\|} = \frac{\sigma\sqrt{n}}{\|\mB\|}.
\end{equation}
The above quantity (rather than the deviation $\sigma$) is to measure the model's noise level.

\begin{remark}
Our goal is to simply perform the symmetric Tucker decomposition of the sample moment tensor $\tMs \in \S^d(\Real^n).$ This is different from inference procedures on the linear factor model, e.g., independent component analysis (ICA) \cite{LMV002}. 
\end{remark}

\begin{remark}
Besides the model \eqref{eqn: factor model} or other heavy-tailed datasets, we emphasize that our proposed algorithms are able to efficiently decompose $\tMs$ when the sample vectors come from any probabilistic model.
\end{remark}

All numerical tests are implemented using MATLAB R$2022$a on a MacBook Pro with 16GB RAM memory. 
Each reported plot is averaged over $5$ independent simulations with independent, random initializations.
The matrix $\mB$ in  \eqref{eqn: factor model} is randomly generated with i.i.d. standard normal entries and then fixed for all runs. 
 We generate new $\vf$ and $\pmb{\eps}$ to form the data samples in each stochastic simulation.
 The parameters of each experiment, namely $n$, $r_{\text{true}}$, $r$, $p$, $d$, $\SNR^{-1}$, are specified below.

\subsection{ Implicit versus explicit algorithms}\label{subsec: implici vs explicit}

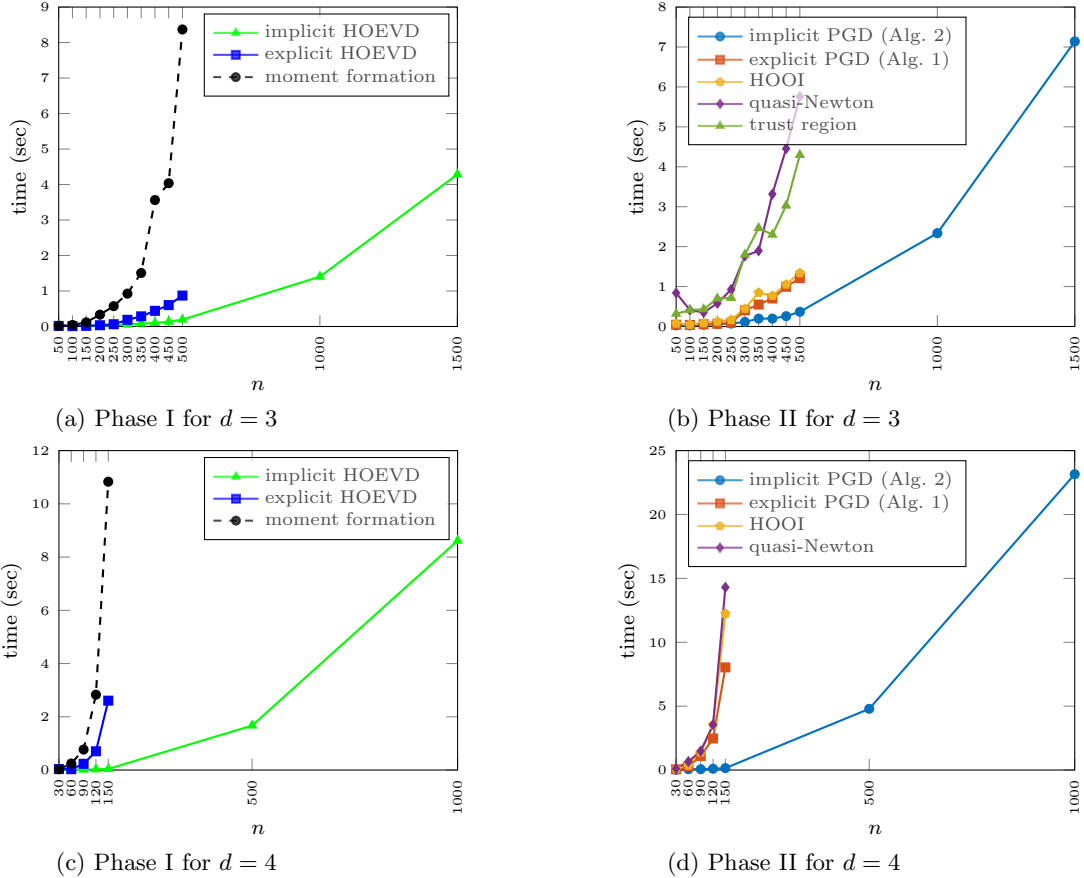
\begin{figure}[!ht]
  \definecolor{c1}{HTML}{0072BD}%
  \definecolor{c2}{HTML}{D95319}%
  \definecolor{c3}{HTML}{EDB120}%
  \definecolor{c4}{HTML}{7E2F8E}%
  \definecolor{c5}{HTML}{77AC30}%
  \pgfplotsset{
    axis style={
      width=5.3cm,height=4.25cm,
      xlabel=$n$,
      ylabel={time (sec)},
      plot style={thick,mark size=1.5pt},
      font=\footnotesize,
      ymin=0,
      label style={inner sep=0em},
      ticklabel style={font=\tiny},
      xticklabel style={rotate=90},
      legend style={font=\scriptsize,fill opacity=0.7,text opacity=1,row sep=-2pt,legend cell align=left},
      xtick=data,
      numformat={1000 sep=},
    },
    mf/.style={black,dashed,mark=*},
    ehoevd/.style={blue,mark=square*},
    ihoevd/.style={green,mark=triangle*},
    ipgd/.style={c1,mark=*},
    epgd/.style={c2,mark=square*},
    hooi/.style={c3,mark=pentagon*},
    qn/.style={c4,mark=diamond*},
    tr/.style={c5,mark=triangle*},    
  }%  
  \centering
  \begin{subfigure}{0.49\textwidth}\flushright
    \begin{tikzpicture}
      \begin{axis}[
          xmin=50,xmax=1500,ymax=9,ytick distance=1,
          legend entries={implicit HOEVD, explicit HOEVD, moment formation}
        ]
        \addplot[ihoevd] table[x=n,y=implicit HOEVD time (avg)]{fig2_a.csv};
        \addplot[ehoevd] table[x=n,y=explicit HOEVD time (avg)]{fig2_a.csv};
        \addplot[mf] table[x=n,y=MoM construction time]{fig2_a.csv};
      \end{axis}
    \end{tikzpicture}
    \caption{Phase I for $d=3$}
    \label{fig:plot2a}
  \end{subfigure}
  \begin{subfigure}{0.49\textwidth}\flushright
    \begin{tikzpicture}
      \begin{axis}[
          xmin=50,xmax=1500,ymax=8,
          ytick distance=1,
          legend entries={implicit PGD (Alg.~2), explicit PGD (Alg.~1),
            HOOI, quasi-Newton, trust region},
          legend pos=north west,
        ]
        \addplot[ipgd] table[x=n,y=implicit PGD time (avg)]{fig2_b.csv};
        \addplot[epgd] table[x=n,y=explicit PGD time (avg)]{fig2_b.csv};
        \addplot[hooi] table[x=n,y=HOOI time (avg)]{fig2_b.csv};
        \addplot[qn] table[x=n,y=QN time (avg)]{fig2_b.csv};
        \addplot[tr] table[x=n,y=trust region time (avg)]{fig2_b.csv};
      \end{axis}
    \end{tikzpicture}
    \caption{Phase II for $d=3$}
    \label{fig:plot2b}
  \end{subfigure}
  
  \begin{subfigure}{0.49\textwidth}\flushright
    \begin{tikzpicture}
      \begin{axis}[
          xmin=30,xmax=1000,ymax=12,
          ytick distance=2,
          legend entries={implicit HOEVD, explicit HOEVD, moment formation}
        ]
        \addplot[ihoevd] table[x=n,y=implicit HOEVD time (avg)]{fig2_c.csv};
        \addplot[ehoevd] table[x=n,y=explicit HOEVD time (avg)]{fig2_c.csv};
        \addplot[mf] table[x=n,y=MoM construction time]{fig2_c.csv};
      \end{axis}
    \end{tikzpicture}
    \caption{Phase I for $d=4$}
    \label{fig:plot2c}
  \end{subfigure}
  \begin{subfigure}{0.49\textwidth}\flushright
    \begin{tikzpicture}
      \begin{axis}[
          xmin=30,xmax=1000,ymax=25,
          ytick distance=5,
          legend entries={implicit PGD (Alg.~2), explicit PGD (Alg.~1),
            HOOI, quasi-Newton, trust region},
          legend pos=north west,
        ]
        \addplot[ipgd] table[x=n,y=implicit PGD time (avg)]{fig2_d.csv};
        \addplot[epgd] table[x=n,y=explicit PGD time (avg)]{fig2_d.csv};
        \addplot[hooi] table[x=n,y=HOOI time]{fig2_d.csv};
        \addplot[qn] table[x=n,y=QN time]{fig2_d.csv};
      \end{axis}
    \end{tikzpicture}
    \caption{Phase II for $d=4$}
    \label{fig:plot2d}
  \end{subfigure}  
  \caption{The scalability of implicit and explicit algorithms. Implicit method is much
    more scalable.}
  \label{imp_exp_time}
\end{figure}

We illustrate the scalability of implicit and explicit implementations in \cref{imp_exp_time}.
The experiments are \emph{non-streaming}. 
The implicit methods run on a fixed set of observations $\mX \in \Real^{n \times p}$.
The explicit approaches require constructing the sample moment tensor $\tMs\in\S^d(\Real^n)$ as overhead, and then operating on it. 
Due to resource limitations, the maximal tensor dimension $n$ for explicit methods is smaller than maximal dimension for the implicit methods. We set the sample size to be $p = 5n$ and choose the rank as $r = r_{\text{true}} = 5.$
We consider non-streaming methods, so $b_1 = b_2 = p$.
The sample moment formation uses the \texttt{symktensor} package in the Tensor Toolbox for MATLAB~\cite{tensortoolbox}.
The noise level is $\SNR^{-1} = 0.5,$ see \eqref{eqn: SNR}.

%%% Local Variables:
%%% mode: latex
%%% TeX-master: "fig2_standalone"
%%% End:

We conduct the experiments in two phases for order-$3$ and $4$ sample moments $\tMs$. 
For Phase I (initialization), we compare the efficiency of the implicit HOEVD result,
\begin{equation}
    \label{eqn: implicit hoevd}
\mQ_{\text{hoevd}}\in \st(n,r) ~\text{given~by~the~leading}~ r ~\text{eigenvectors~ of}~ \mX \left(\mX^\top \mX\right)^{[d-1]}\mX^\top\!/p^2,
\end{equation}
and its explicit counterpart \eqref{eqn: hoevd}. The equivalence of the two methods is by \eqref{eqn: implicit 2}.

Then in Phase II (iterations), we run the SPGD \cref{algorithm} and various explicit algorithms, namely the explicit PGD \cref{basic PGD}, HOOI \cite{LMV001}, quasi-Newton \cite{SL10} and trust-region \cite{IAHL11} methods.\footnote{Similar to our algorithm scheme, the comparable algorithms are mostly gradient and Hessian-driven methods. Since we could not identify an implementation for the Jacobi rotation method \cite{IAV13,LUC18}, it is not included in the comparison.}

\sloppy 
For the last three algorithms, we use the code \texttt{tucker\_als.m} in Tensor Toolbox \cite{tensortoolbox}, \texttt{symalgQNlc.m} in package \cite{S09} and \texttt{lmlra3\_rtr.m} in TensorLab \cite{tensorlab}, respectively. 
HOOI and trust-region apply to general asymmetric tensors.
We omit the trust-region method in $d=4$ case as its code does not apply to the fourth order. 
The step size for \cref{algorithm} is $c_2=1000$.
Note that the implicit \cref{algorithm} matches the explicit \cref{basic PGD} iteration-by-iteration, because there is no streaming.  However, they calculate the updates differently because \cref{algorithm} is implicit.

\cref{imp_exp_time} shows the running time of the two phases as the dimension $n$ increases. 
The time on the vertical axes shows when the algorithm reaches a relative gradient to be less than $10^{-12}.$ 
This measure of convergence precision is used in \cite{ES09,SL10,IAHL11} and defined as:
\[
\texttt{relative~gradient} := \frac{\|\text{grad}~F(\mQ)\|}{F(\mQ)}.
\]
Here
$
\grad F(\mQ) = (\id_n - \mQ\mQ^\top)\nabla F(\mQ) \in \Real^{n \times r}
$
denotes the Riemannian gradient on the Stiefel manifold \eqref{eqn: stiefel}, where $\nabla F(\mQ)\in \Real^{n \times r}$ is the Euclidean gradient  \eqref{eqn: euclid grad}. 
We use the Riemannian gradient $\grad F(\mQ)$ because the possible update direction on the Stiefel manifold is constrained to the tangent  
space at $\mQ.$
The time for checking the relative gradient is included in the runtime.

 It is evident from \cref{imp_exp_time} that the implicit HOEVD \eqref{eqn: implicit hoevd} and SPGD (\cref{algorithm}) are much faster than the explicit algorithms. 
 The advantage of being tensor-free  allows the implicit implementation to handle computation in higher dimensions. Meanwhile, the prohibitive upfront and storage costs in the explicit methods restrict their scalability.

\subsection{Streaming versus full data usage}
\label{subsec: streaming vs full}
\begin{figure}[!ht]
  \centering
  \definecolor{c1}{HTML}{0072BD}%
  \definecolor{c2}{HTML}{D95319}%
  \definecolor{c3}{HTML}{EDB120}%
  \definecolor{c4}{HTML}{7E2F8E}%
  \definecolor{c5}{HTML}{77AC30}%
  \pgfdeclarelayer{background}
  \pgfdeclarelayer{inbetween}
  \pgfsetlayers{background,inbetween,main}  
  \pgfplotsset{
    axis style={
      width=5.5cm, height=4.25cm, ymin=0,
      xlabel={time (sec)}, ylabel={subspace error},
      plot style={table/x=time (average), table/y=error (average), thick, mark size=1.5pt},
      font=\footnotesize, label style={inner sep=0em}, ticklabel style={font=\tiny,inner sep=0.2em},
      legend style={font=\scriptsize,legend cell align=left}, ynumformat={fixed,precision=1,zerofill},
    },
    every node near coord/.append style={fill=white, fill opacity=0.9, text opacity=1,font=\scriptsize, inner sep=0.1em, yshift=1.2mm,
      /pgf/number format/fixed,/pgf/number format/precision=3,/pgf/number format/zerofill},
    b25/.style={c1,mark=triangle*},
    b50/.style={c2,mark=diamond*},
    b100/.style={c3,mark=*},
    full/.style={c4,mark=square*},
    de/.style={c5,only marks,mark=pentagon*,mark size=2pt,mark style={draw=black,thin}},lastpoint/.style={nodes near coords,unbounded coords=discard,x filter/.expression={\coordindex == #1 ? x : nan}},
    lastpoint/.style={nodes near coords,unbounded coords=discard,
      x filter/.expression={\coordindex == #1 ? x : nan}},
    lastpoint/.default=10,    
    shifty/.style={every node near coord/.append style={#1}}
  }%
  \tikzset{subcaption/.style={text width=6cm,align=flush left,anchor=north,font=\footnotesize,yshift=-3mm,inner ysep=-3mm}}%
  \begin{tikzpicture}
    \begin{pgfonlayer}{inbetween}
      \begin{axis}[name=A,
        xmin=0,ymin=0,ymax=3.5,xtick distance=1,
        legend entries={$b_1=25$,$b_1=50$,$b_1=100$,eigensolver},
        legend to name={mylegend A},
        ]
        \addplot[b25] table {fig3a_b25.csv};
        \addplot[b50] table {fig3a_b50.csv};
        \addplot[b100] table {fig3a_b100.csv};
        \addplot[de,lastpoint=0,shifty={above left,xshift=2mm,name=foo}]
        table {fig3a_eig.csv};
        \coordinate (AAA) at (axis cs:1.3,0);
      \end{axis}      
    \end{pgfonlayer}
    \begin{pgfonlayer}{background}
      \fill[gray!25] (A.north west) rectangle (AAA);
    \end{pgfonlayer}
    \path[fill=white] (A.south west) ++(1.5,0.75) coordinate (AA);
    \begin{axis}[name=AA,width=3.00cm,anchor=south west,at=(AA),
      xlabel=,ylabel=,
      xmin=0,xmax=1.3,ymin=0,ymax=3.5,
      every outer x axis line/.append style={draw=gray},
      every outer y axis line/.append style={draw=gray},
      ticklabel style={gray},
      ]
      \addplot[b25] table {fig3a_b25.csv};
      \addplot[b50] table {fig3a_b50.csv};
      \addplot[b100] table {fig3a_b100.csv};
      \addplot[b25,lastpoint,shifty={yshift=3.3mm}] table {fig3a_b25.csv};
      \addplot[b50,lastpoint,shifty={yshift=1.5mm}] table {fig3a_b50.csv};
      \addplot[b100,lastpoint,shifty={xshift=-1.5mm}] table {fig3a_b100.csv};
    \end{axis}
    \begin{pgfonlayer}{inbetween}
      \fill[white] (AA.outer south west) rectangle (AA.outer north east);
      \path (AA.north) ++(-0.5,0) node[font=\footnotesize,text width=2cm,align=flush center,text=gray,anchor=north] {Zoom of\\gray area};
    \end{pgfonlayer}
    \path (A.north east) node[anchor=north east]{\pgfplotslegendfromname{mylegend A}};
    \draw[gray!50,dashed] (A.north west) -- (AA.north west);
    \draw[gray!50,dashed] (AAA) -- (AA.south east);

    \path(A.south east) ++(0.4,0) coordinate (B root);
    % ------------------------------
    \begin{pgfonlayer}{inbetween}
      \begin{axis}[name=B,anchor=left of south west,at=(B root),
        xmin=0,ymin=0,ymax=0.4,xtick distance=5,
        legend entries={$b_2=25$,$b_2=50$,$b_2=100$,full},
        legend to name={mylegend B},ylabel={},
        ]
        \addplot[b25] table {fig3b_b25.csv};
        \addplot[b50] table {fig3b_b50.csv};
        \addplot[b100] table {fig3b_b100.csv};
        \addplot[full] table {fig3b_full.csv};
        \addplot[full,lastpoint=3,shifty={xshift=-1.5mm}] table {fig3b_full.csv};
        \coordinate (BBB) at (axis cs:0.55,0);
      \end{axis}      
    \end{pgfonlayer}
    \begin{pgfonlayer}{background}
      \fill[gray!25] (B.north west) rectangle (BBB);
    \end{pgfonlayer}
    \path[fill=white] (B.south west) ++(1.5,0.75) coordinate (BB);
    \begin{axis}[name=BB,width=3.00cm,anchor=south west,at=(BB),
      xlabel=,ylabel=,
      xmin=0,xmax=0.55,ymin=0,ymax=0.4,
      every outer x axis line/.append style={draw=gray},
      every outer y axis line/.append style={draw=gray},
      ticklabel style={gray},
      ]
      \addplot[b25] table {fig3b_b25.csv};
      \addplot[b50] table {fig3b_b50.csv};
      \addplot[b100] table {fig3b_b100.csv};
      \addplot[b25,lastpoint,shifty={yshift=1.6mm}] table {fig3b_b25.csv};
      \addplot[b50,lastpoint,shifty={yshift=0.9mm,xshift=1mm}] table {fig3b_b50.csv};
      \addplot[b100,lastpoint,shifty={xshift=-1.4mm}] table {fig3b_b100.csv};
    \end{axis}
    \begin{pgfonlayer}{inbetween}
      \fill[white,fill opacity=0.9] (BB.outer south west) rectangle (BB.outer east);
      \fill[white] (BB.outer west) rectangle (BB.outer north east);
      \path (BB.north) ++(-0.5,0) node[font=\footnotesize,text width=2cm,align=flush center,text=gray,anchor=north] {Zoom of\\gray area};
    \end{pgfonlayer}
    \path (B.north east) node[anchor=north east]{\pgfplotslegendfromname{mylegend B}};
    \draw[gray!50,dashed] (B.north west) -- (BB.north west);
    \draw[gray!50,dashed] (BBB) -- (BB.south east);
    
    \path (A.below south) node[subcaption] {\subcaption{SHOEVD with $T_1=10$ and $c_1=0.05$. The parameter $b_1$ is the streaming blocksize.}};
    \path (B.below south) node[subcaption] {\subcaption{SPGD with $T_2=10$ and $c_2=1$. The parameter $b_2$ is the streaming blocksize.\\
        For full, $c_2=10$ and $b_2=25{,}000$.}};
    %\draw (current bounding box.south west) rectangle (current bounding box.north east);
  \end{tikzpicture}
  \caption{Streaming versus full data usage on randomly-generated third-order moment tensor
    with $n=500$, $p=25,000$, and $\text{SNR}^{-1}=0.5$, using $r=r_{\rm true}$.
    The final error values for all methods are marked on the plot.}
  \label{str_full_time}
\end{figure}
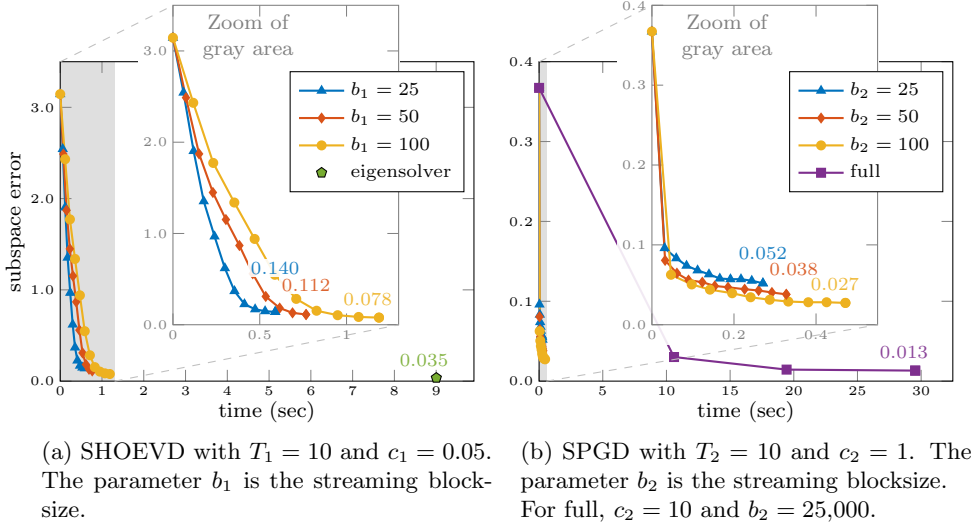

We stress the gains provided by the streaming technique. 
This is demonstrated in stochastic experiments
for (a) SHOEVD (\cref{alg: scalable hoevd}) and (b) SPGD (Phase II of \cref{algorithm})
in \cref{str_full_time}. 
The vertical axes represent the approximation accuracy evaluated by the metric:
\begin{equation*}
\label{eqn: subspace error}
\texttt{subspace~error} := \|\mQ\mQ^\top - \mQ^*\mQ^{*\top} \|.
\end{equation*}
The measure is rotationally invariant. We plot the subspace error against runtime.

%%% Local Variables:
%%% mode: latex
%%% TeX-master: "fig3_standalone"
%%% End:

To compare streaming and full data methods, 
we first generate a  pool of observations $\mX \in \Real^{n \times p}$. 
Each streaming method selects a subset of samples $\mX(:, b(t-1)+1: bt)$ from the pool at step $t$. 
We ensure the total sample number $p$ is just enough for the streaming case with the largest batch size $b_{\max}$. 
The non-streaming counterparts use all the samples $\mX$ at once.
We set the number of iterations $T = 250$, and the total sample size $p = Tb_{\max} = 25,000.$  

In Phase I (a), the full HOEVD refers to the implicit way \eqref{eqn: implicit hoevd},
which is already more efficient than explicit HOEVD \eqref{eqn: hoevd}. 
The step size for all streaming cases of \cref{alg: scalable hoevd} is $c_1 = 0.05.$  
In Phase II (b), all plots start at the same initialization by \cref{alg: scalable hoevd} with $T_1 = 20, b_1 = 50$ and $c_1 = 1.$ 
These parameters for Phase I are different from the (a) setting. They do not affect the efficiency of (b).
The streaming and full cases for the actual Phase II use step size $c_2=1$ and $T_2=10$ respectively.

It is seen in \cref{str_full_time} that the streaming \cref{alg: scalable hoevd,algorithm} have limited loss of accuracy compared to the methods that use the full dataset. 
On the other hand, the runtimes for all streaming cases are nearly negligible compared to timings for the full methods. 

\subsection{Robustness}
\label{subsec: rank and noise}
We test how the noise level in the observations and the choice of target rank affect HOEVD and symmetric Tucker decomposition in \cref{noise,rank}. 
 In the two experiments  below, we plot the evolution of both phases of \cref{algorithm}, see \cref{noise,rank}.   
As the HOEVD phase (in dashed lines) plateaus, its last iterate gives the initialization for  Phase II (in solid lines). 
The number of iterations are shown on the horizontal axes.
The parameters in \cref{algorithm} are set to be $r_{\text{true}} = 5, b_1 = b_2 = 50, c_1 = c_2 = 1.$ 

\begin{figure}[!ht]
  \centering
  \pgfplotsset{
    axis style={
      width=5.15cm, height=4.25cm, xmin=0, ymin=0,    
      xlabel=iteration,
      font=\footnotesize,
      label style={inner sep=0em},
      ticklabel style={font=\tiny,inner sep=0.2em},
      legend style={font=\scriptsize,legend cell align=left,row sep=-2pt},
    },
    table/x expr=\thisrow{iter}-1,
    p1/.style={forget plot,dashed},
    s1/.style={green,mark=triangle*},
    s1a/.style={s1,mark=triangle,p1},
    s2/.style={magenta,mark=diamond*},
    s2a/.style={s2,mark=diamond,p1},
    s3/.style={red,mark=square*,mark size=1.5pt},
    s3a/.style={s3,mark=square,p1},
    s4/.style={blue,mark=*,mark size=1.5pt},
    s4a/.style={s4,mark=o,p1},
  }

  \begin{tikzpicture}
      \begin{axis}[name=A,
          ylabel={subspace error},ytick distance=0.5,
          legend style={at={(0.55,0.95)},anchor=north},
          xtick={40,240},
          xticklabels={$T_1=40$,$T_2=240$}
        ]
        \foreach[count=\cnt] \snr in {0.1 ,0.5 ,1 ,5 }
        {
          \edef\temp{
            \noexpand\addplot[s\cnt a] table[y=si\snr (avg)] {fig4_dash.csv};
            \noexpand\addplot[s\cnt] table[y=si\snr (avg)] {fig4_solid.csv};
            \noexpand\addlegendentry{{$\noexpand\SNR^{-1}$} = \snr}
          }\temp
        }
      \end{axis}
      \path (A.right of south east) ++(0.5,0) coordinate (B root);
      \begin{axis}[name=B,anchor=left of south west,at=(B root),
          ylabel={objective value},
          legend pos=south east,
          ytick distance=1e9,
          xtick={20,120},
          xticklabels={$T_1=20$,$T_2=120$}
        ]
        \addplot[cyan,dashed] expression[domain=0:120] {5.319011E+09};
        \addlegendentry{ground truth}
        \foreach[count=\cnt] \rnk in {3 ,4 ,5 ,6 }
        {
          \edef\temp{
            \noexpand\addplot[s\cnt a] table[y=rank\rnk(avg)] {fig5_dash.csv};
            \noexpand\addplot[s\cnt] table[y=rank\rnk(avg)] {fig5_solid.csv};
            \noexpand\addlegendentry{$r=\rnk$}
          }\temp
        }
      \end{axis}      
     \begin{scope}[every node/.style={text width=5.5cm,align=flush left,anchor=north}]
       \node at (A.below south) {\caption{Two phases of SPGD with different noise levels.
          Tensor has $n=500$, $p=12000$, $d=3$, and  $r=r_{\rm true}=5$.
          In Algorithm 2, $b_1=b_2=50$ and $c_1=c_2=1$.}\label{noise}};      
       \node at (B.below south) {\caption{Two phases of SPGD with different target ranks.
           Tensor has $n=500$, $p=6000$, $d=3$,  and $r_{\rm true}=5$, and $\SNR^{-1}=0.5$.
           In Algorithm 2, $b_1=b_2=50$ and $c_1=c_2=1$.}\label{rank}};      
    \end{scope}
  \end{tikzpicture}  
\end{figure}

%%% Local Variables:
%%% mode: latex
%%% TeX-master: "fig45_standalone"
%%% End:

We first consider the algorithms' sensitivity to noise in the data samples.
To test this, we generate a sequence of clean samples and then plant Gaussian noise with different levels of $\SNR^{-1}$. 
The algorithms' results are shown in \cref{noise}. 
In particular, note the improvement in Phase II over the HOEVD phase. 
The discrepancy is enlarged in the high noise regime ($\SNR^{-1} \geq 1$), when the HOEVD solver provides a poor solution. 
In contrast, relatively good approximation is obtained using the symmetric Tucker decomposition (Phase II) when $\SNR^{-1}$ is no bigger than $1$.

Next, we test various rank choices in \cref{rank}. 
All cases pass over the same set of data observations. 
The algorithms' recovery qualities are reported in terms of the objective value $F$ \eqref{eqn: sym Tucker}, and compared to the ``ground-truth" value $F(\mQ^*)$.
As there is no deterministic moment $\tMs$ in the streaming case, we estimate $F$ using an additional set of testing samples $\mX_{\text{test}} \in \Real^{n \times p_{\text{test}}}$ drawn from  \eqref{eqn: factor model}.  We calculate $F$  implicitly using \eqref{eqn: implicit 1}:
\begin{equation*}
\label{eqn: F}
\begin{array}{ll}
\texttt{objective~value:~} F(\mQ) 
& \displaystyle=\frac{1}{p^2_{\text{test}}}\left\|\sum_{i=1}^{p_{\text{test}}}\left(\mQ^\top \mX_{\text{test}}\right)_i ^{\otimes d} \right\|^2 \\
&= \displaystyle \frac{1}{p^2_{\text{test}}} \left\| \left(\mX_{\text{test}}^\top \mQ \mQ^\top \mX_{\text{test}}\right)^{[d]} \right\|.
\end{array}
\end{equation*}

When $r \geq r_{\text{true}},$ the solutions obtained by \cref{algorithm} fully approximate the moment tensor in the sense that their objective values match the ground truth. 
When $r < r_{\text{true}}$, we can only partially recover the tensor. 
The less the rank is, the lower the objective value. 
Also note that as $r$ decreases, HOEVD's initialization moves further away from the ground truth and the gap between the phases is more pronounced. 

\cref{noise,rank} show that the robustness of the symmetric Tucker decomposition is stronger than that of the HOEVD approximation.
For ``hard" problem instances where  noise is high or the rank is under-specified, the HOEVD phase (in dashed lines) is clearly improved upon during Phase II (in solid lines).

\subsection{Tucker versus CP decomposition}
\label{subsec: Tucker vs CP}
This subsection shows the situations where symmetric Tucker decomposition is more appropriate to use than the symmetric CP counterpart. We provide numerical evidence in the moment tensor approximation quality for low-rank factor model, while maintaining a low storage budget.

We employ the relative error as the accuracy metric: 
\begin{equation*}
\texttt{relative~error} := \frac{\left\|\tMs - \tMs_{\text{test}}\right\|^2  }{\left\|\tMs_{\text{test}}\right\|^2} \times 100 \%.
\end{equation*}
The testing sample moment $\tMs_{\text{test}} \in \S^d (\Real^n)$ serves as the reference solution and is built from the dataset $\mX_{\text{test}} \in \Real^{n \times p_{\text{test}}}$ via \eqref{eqn: sample moment}.
The approximation tensor $\tMs\in \S^d(\Real^n)$ is obtained via either the symmetric Tucker or the symmetric CP \cite{SK20} decompositions. 
These solutions admit the formats:
\begin{equation}
\label{eqn: Tucker and CP}
\tMs_{\text{Tucker}} 
\displaystyle = \tMs_{\text{test}}\cdot\left( \mQ\mQ^\top , \dots, \mQ\mQ^\top \right)\in \S^d (\Real^n) ,\quad
\tMs_{\text{CP}} \displaystyle = \sum_{i=1}^{r'} \lambda_i \v_i^{\otimes d}\in \S^d (\Real^n).
\end{equation}
In \eqref{eqn: Tucker and CP}, we denote the Tucker and CP rank by $r$ and $r'$ respectively. 
The Tucker approximation is determined by the basis $\mQ \in \Real^{n \times r}$ and its core $\tC \in \S^d(\Real^r)$ is via \eqref{eqn: core}.
The CP format consists of the coefficients $\{\lambda_i\}_{i=1}^{r'}\subset \Real $ and the vector components $\{\v_i\}_{i=1}^{r'} \subset \Real^n.$
For the Tucker and CP solutions, we respectively apply implicit stochastic methods: our \cref{algorithm} and Algorithm \texttt{cp\_isym.m} \cite{SK20} with Adam optimization from Tensor Toolbox \cite{tensortoolbox}.

The data observations are drawn from the linear factor model \eqref{eqn: factor model}. Here, the latent factor $\vf \in \Real^{r_{\text{true}}}$ follows the standard Gaussian distribution with zero mean and identity covariance. The noise level is $\SNR^{-1} = 0.05.$
In one simulation, we first draw the full data set $\mX_{\text{test}} \in \Real^{n \times p_{\text{test}}}$ with dimension $n=500$ and sample size $p_{\text{test}} = 10000$.  We form order $d=4$ moment tensors.
At each iteration of both algorithms, we compute the gradient based on a sample batch $\mX_t \in \Real^{n \times b}$ with batch size $b = 50$.
Their objectives are evaluated based on the full testing tensor $\tMs_{\text{test}}.$ 
The algorithms take a complete pass over $\mX_{\text{test}}$ in one simulation (for Tucker) and one epoch (for CP). All epochs of one CP simulation operate on the same data $\mX_{\text{test}}$.

In \cref{table: err 4}, we compare the accuracies of the symmetric Tucker versus CP decompositions, in reference to $\tMs_{\text{test}}$.
Regarding the target rank, the Tucker method sets it equal to the true rank: $r_{\text{true}}$. There are two choices for the CP method: $r_{\text{true}}$ and $r_{\text{true}}^{d-1}.$
For each table entry, we report the lowest error from 5 independent simulations. 

\begin{table}[!htp]
\centering
\small
\begin{tabular}{|c|c|c|c|c|c|}
\hline 
Tucker rank: $r = r_{\text{true}}$ &3 & 4 & 5 & 6 & 7\\ 
\hline
relative error ($\%$) & 0.0027  &0.0035  & 0.0039  & 0.0048 & 0.0059 \\ 
\hline
\hline
CP rank: $r' = r_{\text{true}}$ &3 & 4 & 5 & 6 & 7\\ 
\hline
relative error ($\%$)& 32.74 & 42.36 & 45.94  & 52.77 & 57.59  \\ 
\hline\hline
CP rank: $r' = r_{\text{true}}^{d-1}$ & 27 & 64 & 125 & 216 & 343\\ 
\hline
relative error ($\%$) &0.0025  & 0.0060  &0.0084    & 0.0089 & 0.0106  \\ 
\hline 
\end{tabular}
\caption{\small The error comparisons of the symmetric Tucker and CP decompositions when $n = 500, d = 4.$ The row indicates the target rank choice.
The moment tensor estimation errors are displayed along the columns.
 In the symmetric Tucker method (\cref{algorithm}), Phase I and II have $T_1 = 20, T_2 = 180$ iterations. AdaGrad step sizes are $c_1 = c_2 = 1.$ In the CP algorithm, the iteration number per epoch is $200.$ The maximal epoch number is $100.$  Adam learning rate is $0.05$ and the decay rate is $0.1.$}
\label{table: err 4}
\end{table}

We see that the Tucker estimation achieves the greatest accuracy on a low rank choice.  On the other hand, CP decomposition suffers substantially larger error (10000 times worse than Tucker) when its target rank is small. CP manages to fit moment tensors with high precision similar to Tucker's level for $r' = r_{\text{true}}^{d-1}$. 
However, comparing the storage costs, Tucker is much more economical in representing the moment tensor.
From \eqref{eqn: Tucker and CP}, Tucker costs $rn + \binom{r+d-1}{d} \approx r_{\text{true}} n$ to store $\mQ$ and $\tC$, whereas the CP decomposition costs $r'(n+1) \approx r_{\text{true}}^{d-1} n $ in storage to achieve a similar approximation.

We conclude that there are important situations where the Tucker structure provides a more natural and more accurate format for moment tensors than CP does.  This holds especially when the data samples lie close to a low-rank subspace.

%%% Local Variables:
%%% mode: latex
%%% TeX-master: "../manuscript"
%%% End:

\section{Applications}
\label{sec: applications}
This section illustrates the  applicability of moment Tucker decompositions to real datasets from anomaly detection and portfolio allocation.

The computations are not wholly streaming, as we need to store and preprocess the data first. 
However, only data batches are used while running \cref{alg: scalable hoevd,algorithm}. The task results are based on one simulation. 

\subsection{Anomaly detection in hyperspectral imaging}
\label{subsec: anomaly detection}
Outliers can be discovered using higher-order statistics, when the mean and covariance information are not sufficient \cite{DGP18,A19}. 
In hyperspectral imaging (HSI) \cite{P09}, the detection goal is  to find a low-rank feature subspace, denoted by $\mQ \in \st(n,r)$.  
After transforming HSI pixels by $\mQ,$ unusual targets in the image are revealed \cite{GSJZ14}. 
Intuitively if the anomalies are captured by $\colspan(\mQ)$, they should stand out as bright spots in the transformed image.

This section tests three methods to extract the feature subspace from HSI data:
 Tucker factorization of  the skewness (order-3 standardized moment) \cite{M70}, HOEVD applied to the skewness, and  eigendecompostion of the covariance. 
Here, skewness and covariance are formed from HSI pixels.
The results are qualitatively and quantitatively evaluated according to their performance in capturing anomalies, as in \cite{GSJZ14}.

\begin{figure}[!htb]
\centering
\begin{tikzpicture}[thick,fill,scale=0.8]
  \draw[black,shift={(-1.2,1)}] (0,0) rectangle (2,2);
\node[text width=3cm] at (4,1.4) 
    {band $1 \in \Real^{p_1 \times p_2}$}; 
    \node[text width=3cm] at (3.4,2.2) 
    {$\cdots$}; 
    \node[text width=3cm] at (2.8,2.8) 
    {band $n$}; 
    \node[text width=6cm] at (-3.65,2.5) 
    {pixel $\vx_i = [x_{i,1}, \dots, x_{i,n}]^\top \in \Real^n$}; 
    \draw[black,fill=white,fill opacity=1,shift={(-.6,0.5)}] (0,0) rectangle (2,2);
    \draw[black,fill=white,fill opacity=1] (0,0) rectangle (2,2);
\draw [blue,fill=blue,fill opacity=1](0.4,1.5) node {$\bullet$};
\draw [blue,fill=blue,fill opacity=1](-0.2,2) node {$\bullet$};
\draw [blue,fill=blue,fill opacity=1](-0.8,2.5) node {$\bullet$};
\path[blue,draw] (1,1) node{}  -- (-1.4,3) node{};
\end{tikzpicture}
\caption{HSI hypercube.} 
\label{hypercube}
\end{figure}
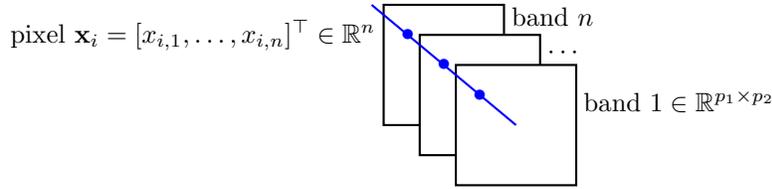

\begin{figure}[!htb]
\centering
\includegraphics[width=.7\textwidth]{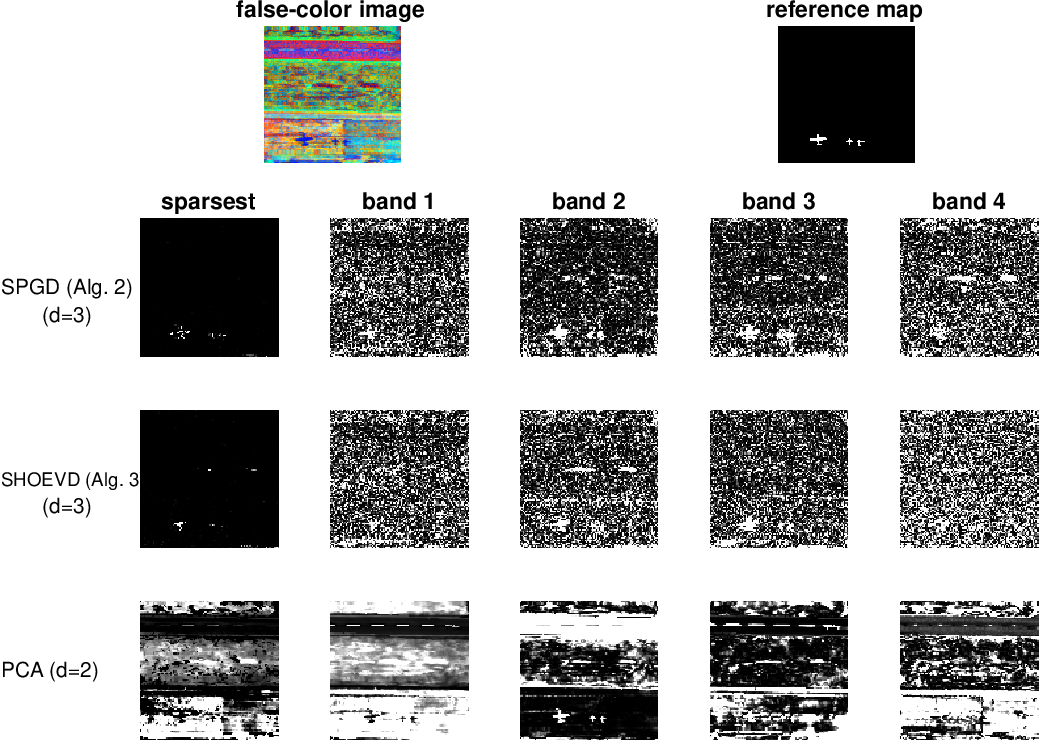}
\caption{{\small Qualitative detection comparison.}}
\label{airport4}
\end{figure}

 \cref{hypercube} shows the data structure of a typical HSI image. It is a $3$-way hypercube $\tX \in \Real^{n \times p_1 \times p_2}.$ The data  consists of $n$ band images in $\Real^{p_1\times p_2}$ and $p = p_1p_2$ HSI pixels in $\Real^n$.
We use a real hyperspectral image from the Airport-Beach-Urban dataset.\footnote{The dataset is available at \url{http://xudongkang.weebly.com/data-sets.html}.} 
The scene was taken at Gulfport, Florida by the Airborne Visible/Infrared Imaging Spectrometer (AVIRIS) sensor. 
The top two images of \cref{airport4} are the false-color composite image and the detection reference image (with white spots depicting the anomalous aircraft). The number of bands is $n=191$ and the size of each band is $p_1\times p_2 = 100 \times 100.$  
We use scalable algorithms for the moment decomposition task.

First, the HSI cube $\tX \in \Real^{ n \times p_1 \times p_2}$ is reshaped into an $n \times p$ matrix.
This is whitened \cite{F87} to give $\mX \in \Real^{n \times p}$.
We choose rank $r=4$, and seek
the leading rank-$r$ subspaces $\mQ \in \Real^{n \times r}$ from the skewness tensor $\tM^3 \in \S^3(\Real^n)$ and the covariance matrix $\mM^2 \in \S^2(\Real^n)$, both formed using the columns in $\mX$.  
The subspaces are found using SPGD (\cref{algorithm}), SHOEVD (\cref{alg: scalable hoevd}), and PCA.  
The parameters for \cref{algorithm,alg: scalable hoevd} are  $b_1 = b_2 = 100$, $c_1 = 0.35, c_2 = 0.5$ and $T_1 = 500, T_2 = 1500.$
For each subspace $\mQ$ obtained, $\mX$ is transformed  into $\mQ^\top\mX \in \Real^{r \times p}$. 
Each row of $\mQ^\top\mX$ is reshaped to a $p_1 \times p_2$ matrix.  
The resulting bands 1 - 4 are shown in \cref{airport4}. 

Since the reference map in \cref{airport4} is sparse, we further pursue a sparsest detection image by optimizing over weights on the bands.
This is achieved by a constrained $\ell_1$-minimization:
\begin{equation*} \label{eq:sparse-ell1}
\vw^{*}=\arg\min_{\vw\in\Real^r}\|\mX^\top\mQ \vw\|_1,~~~~\text{subject~to}~~\|\vw\|_2=1.
\end{equation*}
\noindent  This is solved by the MATLAB Optimization Toolbox function \texttt{fmincon}.  Then we reshape  $\mX^\top\mQ \vw^{*}$ into a band image.  
The resulting sparse images are in the first column of \cref{airport4}.

 Inspecting \cref{airport4}, we see that SPGD and SHOEVD based on skewness distinguish between abnormal and normal pixels well. 
 The subspace provided by SPGD enables especially clear detection, since the targets are the most pronounced. 
It outlines both aircraft individually in the sparse image. 
The covariance information, however, fails to identify anomalies. 

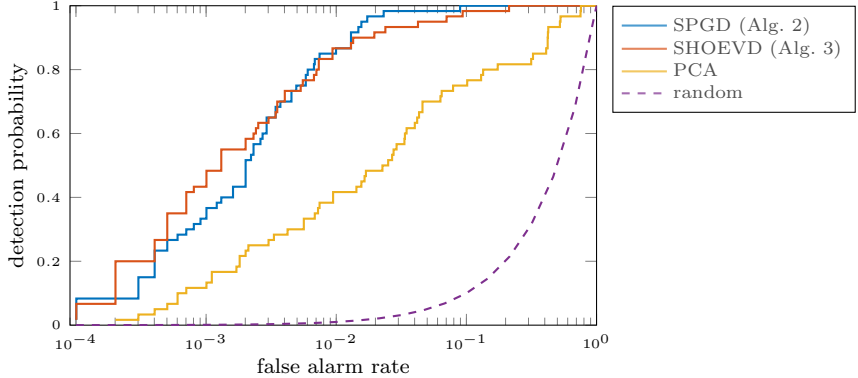
\begin{figure}[!ht]
  \definecolor{c1}{HTML}{0072BD}%
  \definecolor{c2}{HTML}{D95319}%
  \definecolor{c3}{HTML}{EDB120}%
  \definecolor{c4}{HTML}{7E2F8E}%
  \definecolor{c5}{HTML}{77AC30}%
  \pgfplotsset{
    axis style={
      width=7cm,height=4.25cm,
      xlabel=false alarm rate,
      ylabel=detection probability,
      plot style={thick},
      font=\footnotesize,
      label style={inner sep=0em},
      ticklabel style={font=\tiny},
      legend style={font=\scriptsize,fill opacity=0.7,text opacity=1,row sep=-2pt,legend cell align=left},
      xmin=0.9e-4,xmax=1,
      ymin=0,ymax=1,
      legend pos=outer north east,
    },
  }%  
  \centering
  \begin{tikzpicture}
    \begin{semilogxaxis}[
        legend entries={SPGD (Alg.~2), SHOEVD (Alg.~3), PCA, random}
      ]
      \addplot[c1,const plot] table{fig8_alg2_new.csv};
      \addplot[c2,const plot] table{fig8_alg3_new.csv};
      \addplot[c3,const plot] table{fig8_pca_new.csv};
      \addplot[c4,dashed] expression[domain=1e-4:1] {x};
    \end{semilogxaxis}
  \end{tikzpicture}
  \caption{{\small Quantitative detection comparison.}}
  \label{roc}
\end{figure}

 In \cref{roc}, we show a quantitative evaluation in anomaly detection,  via the receiver operating characteristic-area under curve (ROC-AUC).
This quantifies how accurate a method is in classifying normal and abnormal objects, where a higher AUC value indicates better detection. 
To draw the ROC curve, we whiten the matrix $\mX^\top\mQ$ and assign the anomaly score for each pixel as its $\ell_2$-norm. 
The dashed curve indicates a random classifier.
Note that the ROC of skewness-based detection schemes lie strictly above that of PCA. 
The AUC values corresponding to SPGD, SHOEVD and PCA are $0.9947, 0.9900$ and $0.8893$.

%%% Local Variables:
%%% mode: latex
%%% TeX-master: "fig8_standalone"
%%% End:

\subsection{Moment estimation for portfolio allocation}
\label{subsec: portfolio}
Moment estimation is sometimes improved by enforcing a low-rank structure \cite{BCV20}.
In portfolio allocation, we seek a distribution over investments to maximize future returns. 
 Judicious predictions can be obtained by including higher-order moment estimations in allocation optimization, as financial data is  heavy-tailed \cite{JR06}.

Higher moments of financial data should roughly admit a low-rank symmetric Tucker format, because of the Fama-French model \cite{FF93}.  In this section, we test two ways to calculate moments: the sample moment $\tMs$ \eqref{eqn: sample moment} and the low-rank symmetric Tucker approximation to $\tMs$. 
We evaluate the moments by deriving portfolios based on them and comparing the portfolio performances, as in \nolinebreak \cite{BCV20}.

We use the real HFR hedge fund dataset.\footnote{The HFRX indices dataset can be found at: \url{https://www.hfr.com/family-indices/hfrx}.} 
It consists of daily return rates (in percentage) of $n =30$ HFRX indices from 01/02/2009 to 12/31/2019. 
The time series is split into in-sample data (2009 - 2016) $\mX_{\text{in}}\in\Real^{n \times p_{\text{in}}} ~(p_{\text{in}} = 2016)$ and out-of-sample data (2017 - 2019) $\mX_{\text{out}}\in\Real^{n \times p_{\text{out}}} ~(p_{\text{out}} = 754)$. 
The in-sample data is centered. 
We calculate the averaged sample moments \eqref{eqn: sample moment}: $\mM^2 \in \S^2(\Real^n), \tM^3 \in \S^3(\Real^n)$ and $\tM^4 \in \S^4(\Real^n)$. 
Then their low-rank Tucker approximations are computed by SPGD (\cref{algorithm}). 
The rank is set to $r=15$, and  \eqref{eqn: core} is used to solve for the cores once the subspaces $\mQ$ are determined. 
We run Phase II of \cref{algorithm}  twice. 
Starting from a random initialization, we obtain $\mQ$ for $\mM^2$.  This is used an initialization for decomposing $\tM^3$. Finally we tackle $\tM^4$, starting from the $\mQ$ for $\tM^3$. 
The parameters are set to $c_2 = 1, 10^{-2}, 10^{-5}$; $T_2 = 1000, 1000, 200$; and $b_2 = 500, 500, 500$ for the three runs respectively.
This produces $\mM_{\operatorname{Tucker}}^2 \in \S^2(\Real^n), \tM_{\operatorname{Tucker}}^3 \in \S^3(\Real^n)$ and $\tM_{\operatorname{Tucker}}^4 \in \S^4(\Real^n)$.

Given moment estimates, we select the allocation weights that maximize the expected utility (EU) \cite[Supplementary Material, Eqn. (17)]{BCV20}:
\begin{multline}
%\vspace*{-6mm}
\label{eqn: eu}
\begin{array}{ll}
\displaystyle{\vw^*}&:=\displaystyle \arg \max_{\sum_i w_i =1 \atop w_i \geq 0}\EU(\vw)
\, \\
& \equiv \displaystyle  -\frac{\mu}{2} \left\langle\E\left((\vx-\overline{\vx})^{\otimes 2}\right), \vw^{\otimes 2}\right\rangle\displaystyle + \frac{\mu\, (\mu+1)}{6} \left\langle\E\left((\vx-\overline{\vx})^{\otimes 3}\right),\vw^{\otimes 3} \right\rangle  \\
&\displaystyle - \frac{\mu\, (\mu+1)\, (\mu+2)}{24} \left\langle\E\left((\vx-\overline{\vx})^{\otimes 4}\right), \vw^{\otimes 4}\right\rangle,
\end{array}
\end{multline}
where the risk aversion parameter $\mu$ is set to be 1.
We estimate $ \E\left((\vx-\overline{\vx})^{\otimes d}\right) \approx  \tM^d$ in \eqref{eqn: eu},  and then apply the MATLAB Optimization Toolbox function \texttt{fmincon} to output the optimal weights. Meanwhile, estimating $ \E\left((\vx-\overline{\vx})^{\otimes d}\right) \approx \tM_{\operatorname{Tucker}}^d$  results in another weight vector.

\begin{figure}[!ht]
  \definecolor{c1}{HTML}{0072BD}%
  \definecolor{c2}{HTML}{D95319}%
  \definecolor{c3}{HTML}{EDB120}%
  \definecolor{c4}{HTML}{7E2F8E}%
  \definecolor{c5}{HTML}{77AC30}%
  \pgfplotsset{
    axis style={
      every axis title shift=0pt,
      width=10cm,height=5cm,
      xlabel=time,
      ylabel={return rate(\%)},
      title={cumulative daily return (2017--2019)},
      title style={font=\bfseries\footnotesize},
      plot style={thick},
      font=\footnotesize,
      label style={inner sep=0em},
      ticklabel style={font=\tiny},
      legend style={font=\scriptsize,fill opacity=0.7,text opacity=1,row sep=-2pt,legend cell align=left},
      ymin=-4,ymax=8,
      legend pos=south west,
      xtick={0, 250, 502, 754},
      xticklabels={01/03/17, 01/02/18, 01/02/19, 01/02/20},
      enlarge x limits=false,
      ytick distance=2,
      ymajorgrids=true,
    },
  }%  
  \centering
  \begin{tikzpicture}
    \begin{axis}[
      legend entries={low-rank moment, sample moment, equally weighted}
      ]
      \addplot[c1] table{fig9_LowRank.csv};
      \addplot[c2,densely dashed] table{fig9_Sample.csv};
      \addplot[c3,densely dotted] table{fig9_EqWght.csv};
    \end{axis}
  \end{tikzpicture}
  \caption{{\small Portfolio returns under different weights.}}
  \label{return}
\end{figure}
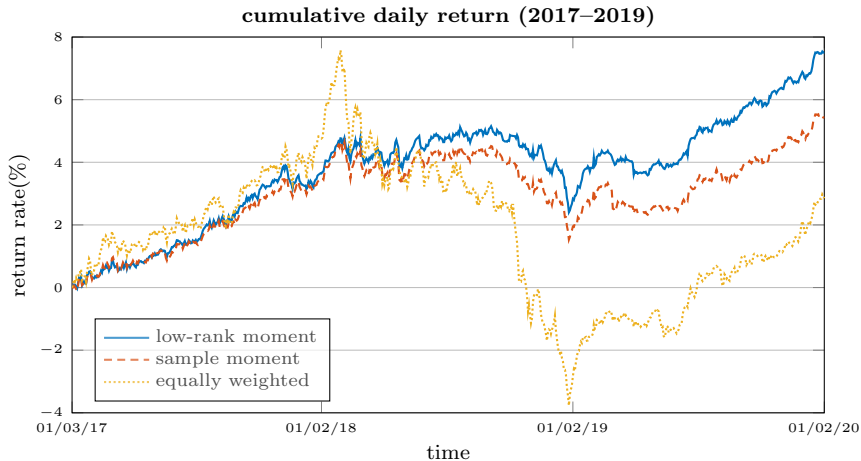

The two portfolios are evaluated on the out-sample data. 
The acquired daily returns $\{y_t \} \subseteq \Real$ are collected in a vector $\vy = \mX_{\text{out}}^\top\vw^* \in \Real^{p_{\text{out}}}.$ 
In \cref{return}, we plot the cumulative daily return rates
$(\prod_{t'=1}^t (1+\frac{y_{t'}}{100})-1)\times 100$
at days $t  = 1, \dots, p_{\text{out}}$.  
We observe that the low-rank moment estimates produce a better return rate than the plain sample moments do. 
Despite mild fluctuations, the two strategies always produce positive return rates. 
A naive strategy, the equally-weighted portfolio, is the most unstable.  It sometimes yields negative returns.

%%% Local Variables:
%%% mode: latex
%%% TeX-master: "fig9_standalone"
%%% End:

Other statistical measures  are shown in \cref{performance}. 
The annual geometric mean of returns over the years $2017 - 2019$ is
$
\left(\left(\prod_{t=1}^{ p_{\text{out}}} (1+\frac{y_{t}}{100})\right)^{1/3}-1\right)\times 100.
$
The centered $d$-th moment ($d=2, 3, 4$) of daily returns $\{y_t\}$ is
$
m^d = \sum_{t=1}^{p_{\text{out}}}\frac{(y_t - \overline{y})^d}{p_{\text{out}}}.
$
Large geometric mean and low absolute values of moments are desirable allocation outcomes. In each row, we highlight the preferred score in bold. They are all achieved by the portfolio based on low-rank Tucker approximations to the moments. 

\begin{table}[!htb]
%\vspace*{-6mm}
\caption{Statistical performance of portfolios.}
\begin{center}
\footnotesize{
\begin{tabular}{cccc}
  \toprule[1pt]
& \bf Equally weighted & \bf Sample & \bf Low rank ($r = 15$) \\ [1mm]
annual geometric mean & $9.4143\times 10^{-1}$  & $2.2434$ & $ \bf 2.3268$ \\ 
$m^2$& $2.6914\times 10^{-2}$ & $8.9233 \times 10^{-3}$ & $\bf 8.5933 \times 10^{-3}$ \\ 
$m^3$& $-2.9724 \times 10^{-3}$ &$ -7.7193 \times 10^{-4}$ &  $\bf -6.0692 \times 10^{-4}$\\ 
$m^4$& $3.4507 \times 10^{-3}$& $5.4339 \times 10^{-4}$ &  $\bf 4.7906\times 10^{-4}$\\ 
 \bottomrule[1pt]
\end{tabular}}
\end{center}
\label{performance}
\end{table}

%%% Local Variables:
%%% mode: latex
%%% TeX-master: "../manuscript"
%%% End:

\section{Convergence analysis}
\label{sec: analysis}

This section presents a convergence theory for the basic PGD method (\cref{basic PGD}) applied to symmetric Tucker decomposition \eqref{eqn: sym Tucker}.  
We note several points in the analysis:  1) The tensor $\tX$ needs not arise as a moment tensor; even when it does, we do not analyze the effects of streaming.
2) Proving convergence of the non-streaming algorithm to critical points is needed first, because \eqref{eqn: sym Tucker} is non-convex, with a complex landscape when $d \geq 3$ (recall \cref{rem:eckart-young}). 
3) Our results do not immediately follow from known general properties of plain GD or projected GD. 
We use the specific properties of the Tucker objective function and QR retraction.  

\subsection{Main statements}
Recall the set-up in \cref{subsec:prob-state}.  
Because the Tucker factorization cost \eqref{eqn: sym Tucker} is non-unique as 
\begin{equation*}
   \left\| \tX \cdot (\mQ, \ldots, \mQ)\right\|^2 = \left\|\tX \cdot (\mQ \mO, \ldots, \mQ \mO)\right\|^2~~\text{for all}~r \times r ~\text{orthogonal matrices}~\mO,
\end{equation*}
the convergence of $\left\{\mQ_t\right\}$ is not relevant.
Instead we focus on the convergence of the projectors, 
$
\mP_t := \mQ_t \mQ_t^\top \in \gr(n,r).
$
Thus, while \cref{basic PGD} runs on the Stiefel manifold and generates the sequence of \emph{bases} $\{\mQ_t\} \subseteq \st(n,r)$, we analyze the sequence of \emph{subspaces} $\{\mP_t\} \subseteq \gr(n,r)$ on the Grassmannian manifold.  

Given scalar step sizes $\{\gamma_t\} \subseteq \Real_{\geq 0}$, the update scheme starts from the PGD iteration on $\mQ$:
\begin{equation}
\label{eqn: Q update}
\mQ_{t+1} = \Phi(\mQ_t) := \qr\left(\mQ_t+ \gamma_t \nabla F(\mQ_t)\right) \in \st(n,r).
\end{equation}
The corresponding update on $\mP$ is defined as 
\begin{equation}
\label{eqn: P update}
\mP_{t+1} = \phi(\mP_t) := \mQ_{t+1}\mQ_{t+1}^\top=\Phi(\mQ_t)\Phi(\mQ_t)^\top \in \gr(n,r).    
\end{equation}

\begin{figure}[!htb]
\centering
\begin{tikzpicture}[thick,fill,scale=0.8]
  \matrix (m) [matrix of math nodes,row sep=2.5em,column sep=2.5em,minimum width=2.5em]
  {
      \mQ_{t+1} = \Phi(\mQ_t) \atop \displaystyle   \st(n,r) & \mP_{t+1} = \mQ_{t+1}\mQ_{t+1}^\top = \phi(\mP_t) \atop \displaystyle  \gr(n,r) \\
     \displaystyle \st(n,r) \atop \mQ_t  & \displaystyle \gr(n,r) \atop  \mP_t = \mQ_t\mQ_t^\top \\};
  \path[-stealth]
    (m-2-1) edge node [left] {$\Phi$} (m-1-1)
        (m-1-1) edge node [below] {} (m-1-2)
            (m-2-2) edge node [right] {$\phi$} (m-1-2)
    (m-2-1) edge node [below] {} (m-2-2);
\end{tikzpicture}
\caption{Update diagram.}
\label{QP update}
\end{figure}

The sequence $\left\{\mP_t\right\}$ is assessed in terms of the cost function
\begin{equation}
\label{eqn:  cost 2}
f(\mP) := \|\tX \cdot\left( \mP, \dots, \mP\right)\!\|^2.
\end{equation}
As noted in \cref{subsec:prob-state}, it is equivalent to the cost in \eqref{eqn: sym Tucker}, as $f(\mP) = F(\mQ)$ when $\mP = \mQ \mQ^{\top}$.
We prove that the iterates $\{\mP_t\}$ converge to a critical point of $f$ on the Grassmannian. 
This is stronger than just showing convergence of the objective values $\{f(\mP_t)\}$ in $\Real$.

The first and second-order criticality conditions are stated as follows in the manifold setting. 
\begin{definition}
A point $\mP \in \gr(n,r)$ is first-order critical for the cost function $f$ \eqref{eqn: cost 2} if
\begin{equation} 
\label{eqn: 1st critical}
\grad f(\mP) = {\bf 0}.
\end{equation}
Here, $\grad f(\mP) \in \tangent_{\mP} \gr(n,r) \subseteq \S^2(\Real^n)$ denotes the Riemannian gradient, which is defined as the Euclidean gradient $ \nabla f(\mP) \in \S^2(\Real^n)$ projected to $\tangent_\mP \gr(n,r)$.
The tangent space at $\mP$ is given \nolinebreak by
\begin{equation}
\label{eqn: tangent space}
\tangent_\mP \gr(n,r) := \left\{ \Delta \mP \in \S^2(\Real^n) ~\vert ~ \mP\Delta \mP+ (\Delta \mP) \mP = \Delta \mP, \, \tr(\Delta \mP) = 0\right\}.
\end{equation}
\end{definition}

\begin{definition}
A point $\mP \in \gr(n,r)$ is second-order critical for the cost function $f$ \eqref{eqn: cost 2} if
\begin{equation}
\label{eqn: 2nd critical}
\grad f(\mP) = {\bf 0}~~~~~~~~\text{and}~~~~~~~~\left\langle \Delta \mP, \hess f(\mP) [\Delta \mP] \right\rangle \leq 0~~\text{for all}~~\Delta \mP \in \tangent_{\mP}\gr(n,r).
\end{equation}
Here, $\hess f(\mP): \tangent_{\mP}\gr(n,r) \to \tangent_{\mP}\gr(n,r)$ is the Riemannian Hessian.  It is calculated in \cite[Eqn.~(2.109)]{HHT07}.
\end{definition}

The first theorem is a uniform guarantee: \cref{basic PGD} always converges to a first-order critical point regardless where the iteration begins.

\begin{theorem}
\label{1st convergence}
Let $\tX \in \S^d(\Real^n)$ be any symmetric tensor, and consider the objective function $f$ defined by \eqref{eqn: cost 2}. 
There exists a constant $\Gamma^*(\tX, r) >0$ such that if the step sizes $\left\{ \gamma_t\right\} \subseteq \Real_{\geq 0}$ satisfy 
$\sup_t \gamma_t \leq \Gamma^*$ and $\inf_t \gamma_t > 0$, the following holds.
For all initializations $\mQ_0 \in \st(n,r)$, the sequence $\left\{ \mP_t \right\} \subseteq \gr(n,r)$ generated via \eqref{eqn: Q update}-\eqref{eqn: P update} by \cref{basic PGD} converges monotonically to a first-order critical point \eqref{eqn: 1st critical} of $f$. 
The convergence is at no less than an algebraic rate.
\end{theorem}

Next, in fact, \cref{basic PGD} converges to a second-order critical point, for almost all initializations.

\begin{theorem}
\label{2nd convergence}
Let $\tX \in \S^d(\Real^n)$ be any symmetric tensor, and consider the cost $f$ defined by \eqref{eqn: cost 2}. 
Use a constant step size $\gamma$ for all iterates.
 There exist a constant $\Gamma^{**}(\tX, r)>0$ and a measure-zero subset $ \mathcal{I} \subseteq \st(n,r)$ such that if $0 < \gamma \leq \Gamma^{**}$, then for all initializations $\mQ_0 \in \st(n,r) \setminus \mathcal{I}$,
 the sequence $\left\{ \mP_t \right\} \subseteq \gr(n,r)$ generated via \eqref{eqn: Q update}-\eqref{eqn: P update} converges to a second-order critical point \eqref{eqn: 2nd critical} of $f$.
\end{theorem}

We build the proofs of \cref{1st convergence,2nd convergence} respectively in  \cref{subsec: 1st convergence,subsec: 2nd convergence}.
The main technical tools are a convergence theorem for real-analytic cost functions (based on the \L ojasiewicz inequality) and the center-stable manifold theorem.  These are stated in \cref{sec:supporting}. 
Detailed proofs for supporting propositions and lemmas are in the Supplementary Materials.  They involve non-trivial calculations for the QR retraction and the symmetric Tucker cost function \eqref{eqn: cost  2}.

\begin{remark} \label{rem:regalia}
The convergence of PGD in $\gr(n,r)$ is not for free.  Previously, Regalia found example tensors $\tX$ where our \cref{algorithm} with a step size approaching infinity shows an oscillatory cycle, rather than convergence.  See \cite[Figures 1 and 4]{R13}. 
\end{remark}

\subsection{Monotonic convergence to first-order critical points}
\label{subsec: 1st convergence}
To establish \cref{1st convergence}, we apply the  guarantee stated in \cref{convergence guarantee}, originally from \cite[Theorem 2.3]{SU15}.  It requires us to prove $\{ \mP_t\}$ meets \textbf{C1} - \textbf{C4} conditions in \cref{convergence guarantee}.

We first compute the Riemannian gradient of $f$ \eqref{eqn: cost 2} on $\gr(n,r)$. 
Although \cref{1st convergence} is in terms of $\mP \in \gr(n,r)$, our  proof uses $\mQ \in \st(n,r)$ due to the explicit update rule in \eqref{eqn: Q update}. 
Thus, we state the relation between $\grad f(\mP)$ and $\mQ, \nabla F(\mQ)$. 
As short-hand notation, define $w: \gr(n,r) \to\S^2(\Real^n)$ by
\begin{equation}
\label{eqn: w}
w(\mP) = d\left\langle\tX\cdot \left(\id_n, \mP, \dots, \mP\right),  \tX\cdot \left(\id_n, \mP, \dots, \mP\right)\right\rangle_{-1} \in \S^2(\Real^n).
\end{equation}
Note that $w(\mP)$ is a positive semidefinite (PSD) matrix.
\begin{proposition}
\label{grad PQ}
For $\mQ\in\st(n,r)$ and $\mP = \mQ\mQ^\top \in \gr(n,r),$ the following holds
\begin{equation}
\label{eqn: grad PQ}
\grad f(\mP) = \sym\Big(2\left(\id_n-\mP\right)w(\mP) \mP \Big)=\sym\Big(\left(\id_n-\mQ\mQ^\top\right)\nabla F(\mQ)\mQ^\top \Big).
\end{equation}
\end{proposition} 
\cref{grad PQ} is proven in \cref{proof grad PQ}. 

Next, we analyze the update rule on $\mP \in \gr(n,r)$. 
For $t\geq 0$ and step size $\gamma_t >0$, 
define $\alpha_t: \gr(n,r) \to \S^2(\Real^n)$ by
\begin{equation}
\label{eqn: alpha}
\alpha_t(\mP) = \mP+ 2\gamma_t \mP w(\mP) \mP \in \S^2(\Real^n).
\end{equation}
 The expression $\alpha_t(\mP)$ and the gradient $\grad f(\mP)$ help to represent the update on $\mP$. 
\begin{lemma} 
\label{difference}
Let a step $t \geq 0$. The difference between the current iterate $\mP_t$ and the next $\mP_{t+1}$ \eqref{eqn: Q update}-\eqref{eqn: P update} is approximated as follows:
\begin{equation}
\label{eqn: difference}
\begin{array}{ll}
\mP_{t+1}-\mP_t &= 4\,\gamma_t \sym\Big(\grad f(\mP_t) \alpha_t(\mP_t)^\dag \Big)\\
&~~+\O\left(\gamma_t^2\left\|\grad f(\mP_t)\alpha_t(\mP_t)^\dag\right\|^2\right)+\O\left(\gamma_t^2\left\|\grad f(\mP_t)\right\|^2\right).
\end{array}
\end{equation}
\end{lemma}
The proof of \cref{difference} is in \cref{proof difference}. It relies on the Taylor expansion of \eqref{eqn: Q update}.

The following lemma gives properties of $\alpha_t(\mP)$ and shows that when the step size is sufficiently small, $\left\|\grad f(\mP_t) \alpha_t(\mP)^\dag \right\|$ and $\left\|\grad f(\mP_t) \right\|$ are on the same order. 
\begin{lemma} 
\label{alpha}
Let $t \geq 0$ and $\mP \in \gr(n,r)$. 
The pseudoinverse $\alpha_t(\mP)^\dag \in \S^2(\Real^n)$ satisfies
$
\alpha_t(\mP)^\dag\alpha_t(\mP) = \mP.
$
Moreover, there exists a constant $\Gamma_1(\tX,r)>0$ such that if $\gamma_t$ satisfies $0< \gamma_t\leq \Gamma_1,$ then $\|\mP - \alpha_t(\mP)^\dag\| \leq \frac{1}{2}$ and 
\begin{equation}
\label{eqn: alpha prop 2}
\frac{\sqrt{2}}{3} \left\|\grad f(\mP_t)\right\| \leq \|\grad f(\mP) \alpha_t(\mP)^\dag \| \leq \frac{\sqrt{2r}}{2}\left\|\grad f(\mP_t)\right\|.
\end{equation}
\end{lemma}
Please refer to \cref{proof alpha} for the proof of \cref{alpha}. It utilizes the fact that $\alpha_t$ is a bounded continuous function of $\mP.$

\begin{proof}[Proof of \cref{1st convergence}]
With the above useful lemmas and proposition in place,
in this proof,
we show that our update scheme for $\mP \in \gr(n,r)$ \eqref{eqn: Q update}-\eqref{eqn: P update} satisfies the conditions \textbf{C1} - \textbf{C4} in \cref{convergence guarantee}. 
Fix a step $t \geq 0$.  Assume $0< \gamma_t \leq \Gamma_1$ where $\Gamma_1$ is as in \cref{alpha}.

\smallskip 

\begin{enumerate}
\item[(\textbf{C1})]
We apply \eqref{eqn: alpha prop 2} that $\O\left(\left\|\grad f(\mP_t)\alpha_t(\mP_t)^\dag\right\|\right)$ has the same scale as 

$\O\left(\left\|\grad f(\mP_t)\right\|\right)$. By \cref{difference}, it follows
\begin{equation}
\label{eqn: difference simple}
\mP_{t+1}-\mP_t = 4\,\gamma_t \sym\Big(\grad f(\mP_t) \alpha_t(\mP_t)^\dag \Big)+\O\left(\gamma_t^2\left\|\grad f(\mP_t)\right\|^2\right).
\end{equation}
It implies that there exists a constant $C_1(\tX, r)>0$  such that
\[
\left\| \mP_{t+1}-\mP_t - 4\,\gamma_t \sym\Big(\grad f(\mP_t) \alpha_t(\mP_t)^\dag \Big)\right\| \leq C_1\gamma_t^2 \left\|\grad f(\mP_t) \right\|^2.
\]
By \eqref{eqn: factor 12} in \cref{proof difference},
\[
\grad f(\mP_t) \alpha_t(\mP_t)^\dag = (\id_n - \mP_t)w(\mP_t) \mP_t \alpha_t(\mP_t)^\dag = (\id_n - \mP_t)w(\mP_t) \alpha_t(\mP_t)^\dag,
\]
and so
\begin{equation*}
\begin{array}{l}
\left\langle \grad f(\mP_t) \alpha_t(\mP_t)^\dag, \alpha_t(\mP_t)^\dag\grad f(\mP_t) \right\rangle \\
=\tr\left( (\id_n - \mP_t)w(\mP_t) \alpha_t(\mP_t)^\dag(\id_n - \mP_t)w(\mP_t) \alpha_t(\mP_t)^\dag  \right)=0.
\end{array}
\end{equation*}
It implies
\[
\left\|\sym\Big(\grad f(\mP_t) \alpha_t(\mP_t)^\dag \Big)\right\| = \displaystyle \frac{\sqrt{2}}{2}\left\|\grad f(\mP_t) \alpha_t(\mP_t)^\dag\right\|. 
\]
By triangle inequality, 
\begin{equation}
\label{eqn: difference 1}
\begin{array}{ll}
\left\| \mP_{t+1}-\mP_t \right\|& \displaystyle \geq 4\times \frac{\sqrt{2}}{2}\gamma_t \left\|\grad f(\mP_t) \alpha_t(\mP_t)^\dag\right\| - C_1\gamma_t^2 \left\|\grad f(\mP_t) \right\|^2\\
& \geq \displaystyle 2\sqrt{2} \times \frac{\sqrt{2}}{3}\gamma_t\left\|\grad f(\mP_t) \right\| - C_1\gamma_t^2 \left\|\grad f(\mP_t) \right\|^2\\
& = \displaystyle\frac{4}{3}\gamma_t\left\|\grad f(\mP_t) \right\| - C_1\gamma_t^2 \left\|\grad f(\mP_t) \right\|^2.
\end{array}
\end{equation}
Here, the second inequality is due to the first half of \eqref{eqn: alpha prop 2}. 
We further constrain the step size as follows
\begin{equation}
\label{eqn: gamma 2}
0< \gamma_t \leq \Gamma_2 :=\min\left(\Gamma_1, \frac{\frac{2}{3}}{C_1\sup_{\mP \in\gr(n,r)}\| \grad f(\mP)\|}\right)\!,
\end{equation}
so that
\begin{equation}
\label{eqn: grad square}
C_1\gamma_t^2 \left\|\grad f(\mP_t) \right\|^2 \leq \frac{2}{3} \gamma_t\left\|\grad f(\mP_t)\right\|.
\end{equation}
Note that $\sup_{\mP \in\gr(n,r)} \| \grad f(\mP) \|$ is finite because $\grad f$ is continuous on $\gr(n,r)$ and the Grassmann is compact.  (We also exclude the trivial case when $\grad f(\mP)\equiv 0$.) Continuing from \eqref{eqn: difference 1} and substituting \eqref{eqn: grad square}, we have
\begin{equation*}
\left\| \mP_{t+1}-\mP_t \right\| \geq (\frac{4}{3}-\frac{2}{3}) \gamma_t \left\|\grad f(\mP_t) \right\| = \frac{2}{3}\gamma_t \left\|\grad f(\mP_t) \right\|.
\end{equation*}
Letting $\kappa = \inf_{t\to \infty}\frac{2}{3}\gamma_t > 0$ ensures that \textbf{C1} \eqref{eqn: C1} holds.  

In preparation for \textbf{C2}, we derive an upper bound for $\left\| \mP_{t+1}-\mP_t \right\|$. The idea is similar to \eqref{eqn: difference 1}:
\begin{equation}
\label{eqn: difference upper bound}
\begin{array}{ll}
\left\| \mP_{t+1}-\mP_t \right\| &\leq 2\sqrt{2}\gamma_t \left\|\grad f(\mP_t) \alpha_t(\mP_t)^\dag \right\| + C_1\gamma_t^2 \left\|\grad f(\mP_t) \right\|^2\\
&\displaystyle \leq 2\sqrt{2} \times \frac{\sqrt{2r}}{2}\gamma_t \left\|\grad f(\mP_t)\right\|+\frac{2}{3}\gamma_t \left\|\grad f(\mP_t) \right\| \\
&\displaystyle\leq \frac{8\sqrt{r}}{3}\gamma_t\left\|\grad f(\mP_t)\right\|.
\end{array}
\end{equation}
The second inequality is obtained from the latter half of \eqref{eqn: alpha prop 2} and \eqref{eqn: grad square}. 

\smallskip

\item[(\textbf{C2})]
Keep the assumption \eqref{eqn: gamma 2} on $\gamma_t$. 
Consider the first-order Taylor's expansion for $f$ on $\gr(n,r)$ \cite[section 10.7]{B22}. 
Inserting \eqref{eqn: difference simple}, the difference of the current and next $f$ is:
\begin{equation}
\label{eqn: f Taylor}
\begin{array}{ll}
f(\mP_{t+1}) -f(\mP_t) & = \langle \grad f(\mP_t), \mP_{t+1} - \mP_t\rangle+ \O\left(\| \mP_{t+1} - \mP_t\|^2\right)\\
& = \left \langle \grad f(\mP_t),4\gamma_t\, \sym\Big(\grad f(\mP_t) \alpha_t(\mP_t)^\dag \Big)\right\rangle \\
&~~+\O\left(\gamma_t^2\max\left(\left\|\grad f(\mP_t) \right\|^3, \left\|\grad f(\mP_t) \right\|^2\right)\right).
\end{array}
\end{equation}

Since $\grad f(\mP_t)$ is symmetric, the first term in the result of \eqref{eqn: f Taylor} is
\begin{equation*}
\begin{array}{l}
\left \langle \grad f(\mP_t),4\gamma_t \sym\Big(\grad f(\mP_t) \alpha_t(\mP_t)^\dag \Big)\right\rangle\\
 = 4\gamma_t \left \langle \grad f(\mP_t),\grad f(\mP_t) \alpha_t(\mP_t)^\dag\right\rangle\\
= 4\gamma_t \left \langle \grad f(\mP_t)\mP_t ,\grad f(\mP_t) \alpha_t(\mP_t)^\dag\mP_t\right\rangle\\
=  4\gamma_t \left \langle \grad f(\mP_t) \alpha_t(\mP_t)^\dag \alpha_t(\mP_t),\grad f(\mP_t) \alpha_t(\mP_t)^\dag\mP_t\right\rangle.
\end{array}
\end{equation*}
The second equation is due to $\alpha_t(\mP_t)^\dag\mP_t = \alpha_t(\mP_t)^\dag.$
The definition \eqref{eqn: alpha} gives $\alpha_t(\mP_t) = \mP_t (\id_n + 2\gamma_t w(\mP_t)) \mP_t$. Continuing from above,
\begin{equation}
\label{eqn: inner prod}
\begin{array}{l}
\left \langle \grad f(\mP_t),4\gamma_t \sym\Big(\grad f(\mP_t) \alpha_t(\mP_t)^\dag \Big)\right\rangle\\
= 4\gamma_t \left \langle \grad f(\mP_t) \alpha_t(\mP_t)^\dag\mP_t (\id_n + 2w(\mP_t)),\grad f(\mP_t) \alpha_t(\mP_t)^\dag\mP_t\right\rangle\\
\geq 4\gamma_t \left\|  \grad f(\mP_t) \alpha_t(\mP_t)^\dag\mP_t \right\|^2 = 4\gamma_t \left\|  \grad f(\mP_t) \alpha_t(\mP_t)^\dag\right\|^2\\
 \geq \displaystyle \frac{8}{9} \gamma_t\| \grad f(\mP_t)\|^2.
\end{array}
\end{equation}
Here we used that $\id_n + 2\gamma_t w(\mP_t) \succeq \id_n$ for $\gamma_t \geq 0$, since $w(\mP_t) \in \S^2(\Real^n)$ is PSD. 
The second inequality is derived from \eqref{eqn: alpha prop 2}. 

Considering the remainder term in \eqref{eqn: f Taylor}, there is a constant $C_2(\tX,r)>0$ such that 
\begin{equation}
\label{eqn: f difference}
\begin{array}{ll}
f(\mP_{t+1}) -f(\mP_t) &\displaystyle \geq\frac{8}{9}\gamma_t\left\| \grad f(\mP_t) \right\|^2\\
& \displaystyle ~~ -C_2\gamma_t^2 \max\left(\left\| \grad f(\mP_t)\right\|^3,\left\|\grad f(\mP_t) \right\|^2\right).
\end{array}
\end{equation}
If the step size $\gamma_t$ is small enough to satisfy
\begin{equation} 
\label{eq:Gamma*}
\begin{array}{ll}
0< \gamma_t \!\!\! & \leq \Gamma^* := \min\left( \Gamma_2, \frac{4}{9C_2}\right),
\end{array}
\end{equation}
where $\Gamma_2>0$ is defined in \eqref{eqn: gamma 2}, then 
\[
 C_2\gamma_t^2 \max\left(\left\|\grad f(\mP_t) \right\|^3,\left\|\grad f(\mP_t)\right\|^2 \right) \leq \frac{4}{9}\gamma_t \left\|( \grad f(\mP_t)\right\|^2.
\]
Inserting into \eqref{eqn: f difference},
\begin{equation*}
\begin{array}{ll}
f(\mP_{t+1}) -f(\mP_t) &\displaystyle\geq  \left(\frac{8}{9}-\frac{4}{9}\right)\gamma_t \left\| \grad f(\mP_t)\right\|^2\\
& =  \displaystyle \frac{1}{6\sqrt{r}} \left(\frac{8\sqrt{r}}{3}\gamma_t \left\|\grad f(\mP_t) \right\|\right) \left\|\grad f(\mP_t) \right\|\\[1em]
& \geq \displaystyle \frac{1}{6\sqrt{r}}\, \| \mP_{t+1} - \mP_t\| \|\grad f(\mP_t)\|.
\end{array}
\end{equation*}
The final inequality is by \eqref{eqn: difference upper bound}. 
Thus \textbf{C2} holds if we set $\sigma =  \frac{1}{6\sqrt{r}}$ in \eqref{eqn: C2}. 

\smallskip

\item[(\textbf{C3})]
``$\Rightarrow$'': From $\grad f(\mP) = \bf 0$ and \eqref{eqn: grad PQ}, we know 
\begin{equation*}
\begin{array}{ll}
 \grad f(\mP)\mQ = {\bf 0} =\grad F(\mQ) = (\id_n-\mQ\mQ^{\top})\nabla F(\mQ),
\end{array}
\end{equation*}
where $\mQ \in \st(n,r)$ is a basis of $\mP \in \gr(n,r)$. It implies that $\colspan(\nabla F(\mQ))$ is in $\colspan(\mQ).$ Note that $\mQ + \gamma\nabla F(\mQ) = \left(\id_n+2\gamma w(\mP)\right)\mQ$ has full column rank, like $\mQ$, because $\id_n+2\gamma w(\mP) \in \S^2(\Real^n)$ is positive definite.  Hence $\colspan\left(\mQ + \gamma\nabla F(\mQ)\right) = \colspan(\mQ)$. Then the next iterate $\Phi(\mQ)$ obtained by QR decomposition still outputs a basis for $\colspan(\mQ).$ So the orthogonal projectors remain the same, i.e. $\phi(\mP) = \mP$.

``$\Leftarrow$'': This directly follows from \textbf{C1} if $\gamma_t$ is sufficiently small to satisfy \eqref{eqn: gamma 2}. 

\smallskip
\item[({\textbf{C4}})]

The Grassmannian is a real-analytic submanifold of $\gr(n,r)$, because $\gr(n,r)$ satisfies property $(3)$ in  \cite[Proposition 2.7.3]{KP02}. Also, the cost function $f$ \eqref{eqn: cost 2} is a polynomial in $\mP \in \gr(n,r)$, so real-analytic as well. By  \cite[Proposition 2.2]{SU15}, all points $\mP\in \gr(n,r)$ satisfy a \L ojasiewicz inequality for $f$.

\end{enumerate}

\smallskip

We have shown that the sequence $\{\mP_t \}$ meets all \textbf{C1 - C4} provided $\{\gamma_t \}$ satisfies  $\inf \gamma_t > 0$ and $\sup \gamma_t \leq \Gamma^*$, with $\Gamma^*$ given in \eqref{eq:Gamma*}. By \cref{convergence guarantee}, it follows \cref{basic PGD} achieves convergence to a first-order critical point of $f$. The proof of \cref{1st convergence} is complete.
\end{proof}

\subsection{Convergence to second-order critical points}
\label{subsec: 2nd convergence}

\sloppy 
We establish \cref{2nd convergence} regarding second-order criticality when the step size is constant.  The proof relies on the center-stable manifold theorem \cite{S87}, which is recalled in \cref{cs}. 
It characterizes the local structure of a smooth local diffeomorphism around a fixed point.

The calculation of the Riemannian Hessian, $\hess f(\mP): \tangent_\mP \gr(n,r) \to \tangent_\mP \gr(n,r)$, and the differential of the update map, $\D \phi(\mP): \tangent_\mP \gr(n,r) \to \tangent_{\phi(\mP)} \gr(n,r)$, are essential to showing the second-order guarantee. 
As short-hand, we define a linear operator $v(\mP):  \tangent_\mP \gr(n,r) \to \Real^{n \times n}$ \nolinebreak by
\begin{equation}
\label{eqn: v}
v(\mP)[\Delta \mP] = d(d-1) \left\langle\tX\cdot \left(\id_n, \Delta \mP, \mP, \dots, \mP\right),  \tX\cdot \left(\id_n, \id_n,\mP \dots, \mP\right)\right\rangle_{-1} \mP \in \Real^{n \times n}.
\end{equation}
\begin{proposition}
\label{hess Dphi formula}
For $\mP \in \gr(n,r)$ and $\Delta \mP \in \tangent_\mP \gr(n,r),$ the Riemannian Hessian of the cost $f$ \eqref{eqn: cost 2} on $\gr(n,r)$ is given by
\begin{equation}
\label{eqn: riem hess f}
\hess f(\mP) [\Delta \mP] =  2\sym\Big((\id_n - \mP) v(\mP)[\Delta \mP] +w(\mP) \left(\Delta \mP \right)\mP-  \Delta \mP w(\mP)\mP\Big).
\end{equation}
The differential of $\phi$ \eqref{eqn: P update} on $\gr(n,r)$ is given by 
\begin{equation}
\label{eqn: D phi}
\begin{array}{l}
\D \phi(\mP) [\Delta \mP] \\
\displaystyle = 2\sym \Big((\id_n -\phi( \mP))\left( \Delta \mP+ 2\gamma w(\mP)\Delta \mP+2\gamma v(\mP)[\Delta \mP] \right) (\mP + 2\gamma w(\mP)\mP)^{\dag}\Big).
\end{array}
\end{equation}
\end{proposition}
The calculations in \cref{hess Dphi formula} are in \cref{proof hess Dphi formula}. These results are derived from standard formula for Riemannian Hessians on $\gr(n,r)$ and the first-order Taylor expansion of $\phi(\mP)$, respectively.
 
 We next relate the eigenvalues of $\hess f(\mP)$ and $\D \phi(\mP)$ when $\mP$ is a first-order critical point of $f$.  In this case, $\mP$ is a fixed point of $\phi$ by {\bf C3} in the proof of \cref{1st convergence}.
 {\sloppy 
\begin{lemma}
\label{hess Dphi}
Let $\mP$ be a first-order critical point \eqref{eqn: 1st critical}.  
If $\hess f(\mP): \tangent_{\mP}\gr(n,r) \to \tangent_{\mP}\gr(n,r)$ has a positive eigenvalue, then $\D \phi(\mP): \tangent_{\mP}\gr(n,r) \to \tangent_{\mP}\gr(n,r)$ has an eigenvalue that is greater than $1$.
\end{lemma}}
\cref{hess Dphi} implies that any first-order, but not second-order, critical point of $f$ is an unstable fixed point of $\phi$. 
The proof is in \cref{proof hess Dphi}. It is needed to express the Jacobian $\D \phi(\mP)$ in terms of an eigenbasis for the Riemannian Hessian $\hess f(\mP).$

The next lemma guarantees that the update map $\phi$ has a nice local property given a sufficiently small step size. 
\begin{lemma}
\label{local diff}
There exists a constant $0 < \Gamma^{**}(\tX,r) \leq \Gamma^*$ (cf. \cref{1st convergence}) such that if the step size satisfies $0 < \gamma \leq \Gamma^{**}$, then the update map $\phi: \gr(n,r) \to \gr(n,r)$ \eqref{eqn: P update} is a local diffeomorphism.  
\end{lemma}

The proof of \cref{local diff} is  in \cref{proof local diff}. We show that $\phi$ is $C^{\infty}$ and $\D\phi(\mP)$ is  invertible for all $\mP\in\gr(n,r)$, when $\gamma$ is small enough.

\begin{proof}[Proof of \cref{2nd convergence}]
Fix the step size $0<\gamma\leq\Gamma^{**}$. From \cref{1st convergence}, the  sequence $\{ \mP_t \}$ always converges to a first-order critical point of $f$ \eqref{eqn: 1st critical}, thus a fixed point of $\phi$ \eqref{eqn: P update} by \textbf{C3} in \cref{subsec: 1st convergence}. In this proof, we consider the case when the limit $\mP \in \gr(n,r)$ is a ``bad'' critical point, meaning it is a first-order critical point but fails the second-order criticality \eqref{eqn: 2nd critical}. 
The idea is to apply \cref{cs} to show that the set of  $\mQ_0$ that converge to bad critical points has zero measure in $\st(n,r)$. 
%We then pull this back to $\st(n,r)$. 

Because $\phi$ is a  local diffeomorphism by \cref{local diff}, we are eligible to apply \cref{cs}.  It gives an open neighborhood $B_\mP \subseteq \gr(n,r)$ around each bad critical point $\mP$ satisfying the conditions in \cref{cs}. 
As $\gr(n,r)$ is second-countable, it has the Lindel\"{o}f property, i.e. every open cover has a countable subcover. So, there exists a countable set $\mathcal{C} \subseteq \gr(n,r)$ consisting of some of the bad critical points such that
\begin{equation}
\label{eqn: countable}
\bigcup_{\mP \in \gr(n,r)  \text{~is~a~bad~critical~point}} B_\mP \quad \quad  = \quad \quad  \bigcup_{\mP\in \mathcal{C}} B_\mP.
\end{equation}

Define the subset
\begin{equation}
\label{eqn: I}
\mathcal{J} \, := \, \bigcup_{\mP \in \mathcal{C}} \bigcup_{t = 0}^\infty \phi^{-t} (W_\mP) \,\,\, \subseteq \,\,\, \gr(n,r).
\end{equation}

We claim that
\begin{equation}
    \left\{ \mP_0 \in \gr(n,r) : \phi^t(\mP_0) \text{~converges~to~a~bad~critical~point~as~} t \rightarrow \infty \right\} \,\, \subseteq \,\, \mathcal{J}.
\end{equation}
To see this, assume $\lim_{t \rightarrow \infty} \phi^t(\mP_0) = \mP$
is bad.
By \eqref{eqn: countable},  $\mP \in B_{\mP'}$ for some $\mP' \in \mathcal{C}$.
Since $B_{\mP'}$ is open, there is $t_0 \geq 0$ such that $\phi^t (\mP_0) \in B_{\mP'}$ for all $t \geq t_0$.
As $\mP'$ is a fixed point of $\phi$, by \eqref{eqn: cs} in \cref{cs}, we obtain $\phi^{t_0} (\mP_0) \in W_{\mP'}$.  
In other words, $\mP_0 \in \phi^{-t_0} (W_{\mP'})$. 
So,  $\mP_0 \in \mathcal{J}$  as claimed.

We claim that $\mathcal{J}$ is measure-zero in $\gr(n,r)$. For each $\mP \in \mathcal{C},$ at least one eigenvalue of $\hess f(\mP)$ is greater than 0. By \cref{hess Dphi}, we know at least one eigenvalue of $\D \phi(\mP)$ is greater than 1. By the first result of \cref{cs}, the dimension of $W_\mP$ is strictly less than the dimension of $\gr(n,r)$.  
So $W_\mP$ is measure-zero in $\gr(n,r)$ by \cite[Corollary 6.12]{L12}. 
Given that the step size satisfies $0 < \gamma \leq \Gamma^{**}$, the map $\phi$ is a local diffeomorphism by \cref{local diff}, and so is $\phi^t$ for all $t \geq 0$. Thus their pre-images of measure-zero sets have zero measure. 
As $(\phi^{t})^{-1} \equiv \phi^{-t}$, the set $\phi^{-t} (W_\mP)$ has zero measure for each $t \geq 0$. 
Then $\mathcal{J}$
is a countable union of measure-zero sets, so has zero measure.

To conclude, 
 define $\mathcal{I} \subseteq \st(n,r)$ as the pre-image of $\mathcal{J} \subseteq \gr(n,r)$ under the map $\mQ \mapsto \mQ \mQ^{\top}$.  
Since the map is a submersion from the Stiefel manifold onto the Grassmannian, $\mathcal{I}$ has zero measure in $\st(n,r)$ as $\mathcal{J}$ does in $\gr(n,r)$.  
For all initializations $\mQ_0$ in  $\st(n,r)$ but outside the null set $\mathcal{I}$, the sequence $\{\phi^t(\mP_0) \} = \{\mP_t \}$ converges to a second-order critical point of $f$ if $0 < \gamma \leq \Gamma^{**}$.
The proof of \cref{2nd convergence} is complete. 
\end{proof}

\section{Conclusion and future work}
\label{sec: conclusions}
In this paper, we developed the PGD method (\cref{basic PGD}) for symmetric Tucker tensor decomposition. It has a simple update scheme compared to other iterative solvers running on Riemannian manifolds. 
We further designed the SPGD method (\cref{algorithm}) and  SHOEVD solver (\cref{alg: scalable hoevd}) to decompose the sample moment tensor \eqref{eqn: sample moment}. The  algorithms are free of constructing and storing moment tensors, and their iterations are conducted by streaming only a small subset of samples. 
Tremendous computational savings were seen in numerical experiments with minimal loss of accuracy.
We also showed the applicability of the Tucker decomposition of moment tensors in comparison to the CP format, as well as to detect anomalies and allocate assets on real datasets.
Finally, utilizing Riemannian manifold optimization theory, we derived theoretical guarantees for the PGD sequence on the Grassmannian.  The sequence always converges to a first-order critical point, and almost surely to a second-order critical point.

There are a number of directions that would be interesting to study in the future.
\begin{enumerate}
\item {\bf Implementation for cumulants.} 
Cumulant tensors are central to non-Gaussian data analysis.
They are nonlinear combinations of moments.
It would be useful to derive scalable implementations for sample cumulants. 
\item {\bf Extending the analysis.} It is certainly worth investigating convergence guarantees in the streaming setting with adaptive step sizes of \cref{algorithm,alg: scalable hoevd}. 
We also want to explore under what conditions is the limiting point guaranteed to be a global optimizer. 

\item {\bf Estimating the core without fresh samples.}  Currently, the core tensor $\tC$ is computed using independent data samples after the subspace $\mQ$ has been calculated.  We would like to figure out a good way of updating $\tC$ while updating $\mQ$, using the same data stream. 
\end{enumerate}

\section*{Acknowledgements}
The authors are especially grateful to Dries Cornilly and Xiurui Geng for their help with the real data. 
We also thank Hemanth Kolla, Amit Singer and Mihai Sirbu for useful discussions. 
R.J. was supported by NSF DMS 1952735. 
J.K. was supported  by start-up grants from the College of Natural Sciences and Oden Institute for Computatational Engineering and Sciences at UT Austin.
R.W. was supported by AFOSR MURI FA9550-19-1-0005, NSF DMS 1952735, NSF HDR-1934932 and NSF 2019844.

%%% Local Variables:
%%% mode: latex
%%% TeX-master: "../manuscript"
%%% End:

\bibliographystyle{siamplain}
\bibliography{references}

\appendix
 \newpage
\section{Supporting theorems} \label{sec:supporting}
\subsection{Convergence guarantee}
\label{subsec: convergence guarantee}
Let $f: \gr(n,r) \to \Real$ be a real-analytic function. 
Consider the optimization  
%\begin{equation*}
$\max_{\mP \in \gr(n,r)} f(\mP).$
%\end{equation*}
Let $\{ \mP_t\} \subseteq \gr(n,r)$ be a sequence. 
We make the following assumptions:
\begin{itemize}
\item
(\textbf{C1}) There exists $\kappa >0$ such that 
\begin{equation}
\label{eqn: C1}
\|\mP_{t+1}-\mP_t\| \geq \kappa\, \|\grad f(\mP_t)\|;
\end{equation}

\item
(\textbf{C2}) There exists $\sigma >0$ such that 
\begin{equation}
\label{eqn: C2}
f(\mP_{t+1}) - f(\mP_t) \geq \sigma \|\mP_{t+1}-\mP_t\|\, \|\grad f(\mP_t)\|;
\end{equation}

\item
(\textbf{C3}) First-order optimality is equivalent to the fixed point condition, i.e. 
\begin{equation}
\grad f(\mP_t) = {\bf 0} \Leftrightarrow \mP_{t+1} = \mP_t;
\end{equation}

\item
(\textbf{C4}) [\L ojasiewicz inequality]
There exist $\delta, \rho >0$ and $\theta \in (0,\frac{1}{2}]$ such that for all first-order critical points $\mP \in \gr(n,r)$ of $f$ \eqref{eqn: 1st critical} and all points $\mP' \in \gr(n,r)$ with $\| \mP'-\mP\|\leq \delta$,  
\begin{equation}
\left\vert f(\mP') - f(\mP) \right\vert^{1-\theta} \leq \, \rho \, \|\grad f(\mP')\|.
\end{equation}
\end{itemize}

\begin{theorem}[Theorem 2.3 in \cite{SU15}]
\label{convergence guarantee}
Under the above set-up, if $\{\mP_t\}$ satisfies conditions \textup{\textbf{C1}} - \textup{\textbf{C4}}, it converges monotonically to a first-order critical point of $f$ at no less than an algebraic rate. 
\end{theorem}

\subsection{Center-stable manifold theorem}
\label{sec: cs}
\begin{theorem}[Theorem III.7(2) in \cite{S87}]
\label{cs}
Let $\phi: \gr(n,r) \to \gr(n,r)$ be a  local diffeomorphism.
Let $\mP \in \gr(n,r)$ be a fixed point of $\phi$. There exist an open neighborhood $B_\mP \subseteq \gr(n,r)$ of $\mP$ and a smoothly embedded disk $W_\mP \subseteq \gr(n,r)$ containing $\mP$ such that
\begin{enumerate} 
\item
the dimension of $W_\mP$ is the number of linearly independent eigenmatrices of $D \phi(\mP)$ with eigenvalues of magnitude no more than 1; 
\item
and 
\begin{equation}
\label{eqn: cs}
\left\{ \mX \in \gr(n,r):  \phi^t(\mX) \in B_\mP,~\forall t \geq 0 \right\}  \, \subseteq \, W_\mP. 
\end{equation}
\end{enumerate}
\end{theorem}

\begin{remark}
Although the  center-stable manifold statement  is discussed for  Euclidean space in \cite{S87}, as noted in \cite[Theorem 1]{LPPSJR19}, it can be extended to any manifold. 
We specify the manifold to be the Grassmannian in this paper.
\end{remark}

In this supplemental document, we present the detailed proofs for all of the supporting propositions and lemmas in \cref{sec: analysis}. 

\section{Proof of \cref{grad PQ}}
 \label{proof grad PQ}
\begin{proof}
We first compute the Riemannian gradient, $\grad f(\mP) \in \tangent_{\mP} \gr(n,r)$. It is the orthogonal projection of the Euclidean gradient $\nabla f(\mP) \in \Real^{n \times r}$ onto the tangent space at $\mP$. 
The projection $\pi_\mP: \S^2(\Real^n) \to \tangent_{\mP} \gr(n,r)$ is computed as follows \cite[Proposition 2.1]{HHT07}:
\begin{equation*}
\pi_\mP(\mX) = \ad^2_\mP(\mX) \in \tangent_{\mP} \gr(n,r)~ \text{with}~ \ad_\mP^2 = \ad_\mP \circ \ad_\mP~\text{and}~ \ad_\mP(\mX) = \mP\mX-\mX\mP
\end{equation*} 
for $\mX \in \S^2(\Real^n)$. Here, $\ad$ stands for the adjoint operator. 

The Euclidean gradient of $f$ \eqref{eqn: cost 2} is 
\begin{equation}
\label{eqn: euclid grad f}
\begin{array}{ll}
\nabla f(\mP) \!\! & = d\left\langle \tX\cdot\left( \id_n, \mP, \dots ,\mP\right),  \tX\cdot\left( \mP, \dots ,\mP\right)\right\rangle_{-1}\\
&~~+ d\left\langle\tX\cdot\left( \mP, \dots ,\mP\right),  \tX\cdot\left( \id_n, \mP, \dots ,\mP\right)  \right\rangle_{-1}\\
& = w(\mP) \mP + \mP w(\mP).
\end{array}
\end{equation}
Hence the Riemannian gradient equals 
\begin{equation*}
\begin{array}{ll}
\grad f(\mP) & = \ad^2_\mP \left(\nabla f(\mP)\right)\\
& = \ad_\mP\left(\mP\nabla f(\mP)-\nabla f(\mP)\mP \right)\\
& = \mP\left(\mP\nabla f(\mP)-\nabla f(\mP)\mP\right) - \left(\mP\nabla f(\mP)-\nabla f(\mP)\mP\right)\mP\\
& = \mP \left(\nabla f(\mP)\right)-2\mP \left(\nabla f(\mP)\right)\mP+\left(\nabla f(\mP)\right)\mP
\end{array}
\end{equation*} 
Substituting \eqref{eqn: euclid grad f} into the above equation,
\begin{equation}
\label{eqn: riem grad fP}
\begin{array}{ll}
\grad f(\mP) & = \mP\left(w(\mP)\mP+\mP w(\mP) \right)\\
& ~~-2\mP\left( w(\mP)\mP+\mP w(\mP) \right)\mP+\left(w(\mP)\mP+\mP w(\mP)\right)\mP\\
& = \mP w(\mP)\left(\id_n-\mP \right) + \left(\id_n-\mP \right) w(\mP) \mP = \sym\Big(2 \left(\id_n-\mP \right) w(\mP) \mP \Big).
\end{array}
\end{equation}

Next, we relate the gradients of $f(\mP)$ and $F(\mQ)$. The Euclidean gradient of $F$ is 
\begin{equation}
\label{eqn: euclid grad FQ}
\begin{array}{ll}
\nabla F(\mQ) &= 2d\left\langle \tX\cdot\left( \id_n, \mQ, \dots ,\mQ\right),  \tX\cdot\left( \mQ, \dots ,\mQ\right)\right\rangle_{-1}\\
& = 2d\left\langle \tX\cdot\left( \id_n, \mP, \dots ,\mQ\right),  \tX\cdot\left( \id_n, \mP, \dots ,\mP\right)\right\rangle_{-1} \mQ\\
& = 2d\left\langle \tX\cdot\left( \id_n, \mP, \dots ,\mP\right),  \tX\cdot\left( \id_n, \mP, \dots ,\mP\right)\right\rangle_{-1}\mQ = 2 w(\mP)\mQ.
\end{array}
\end{equation}
Therefore,
\begin{equation*}
\left(\id_n - \mQ\mQ^\top \right)\nabla F(\mQ)\mQ^\top = 2\left(\id_n - \mP \right) w(\mP)\mQ\mQ^\top =  2\left(\id_n - \mP \right) w(\mP)\mP.
\end{equation*}
By \eqref{eqn: riem grad fP}, we have 
\[
\grad f(\mP) = \sym\Big(\left(\id_n - \mQ\mQ^\top \right)\nabla F(\mQ)\mQ^\top\Big). 
\]

The proof of \cref{grad PQ} is complete. 
\end{proof}

\section{Proof of \cref{difference}}
\label{proof difference}
\begin{proof}
We derive the update on $\mP_t$ \eqref{eqn: P update} based on the iterate $\mQ_t$ \eqref{eqn: Q update}. We decompose the update direction $\gamma_t \nabla F(\mQ_t)$ into two parts respectively in $\colspan(\mQ_t)$ and $\colspan(\mQ_t^\perp)$, i.e.
\[
\gamma_t \nabla F(\mQ_t) = \gamma_t \mQ_t\mQ_t^\top \nabla F(\mQ_t) +\gamma_t\left(\id_n -\mQ_t\mQ_t^\top \right) \nabla F(\mQ_t). 
\]

We approximate $\mQ_{t+1} = \qr\left(\mQ_t +  \gamma_t \nabla F(\mQ_t)\right)$ by the first-order Taylor expansion of $\qr$ centered at
$\mQ_t+ \gamma_t \mQ_t\mQ_t^\top \nabla F(\mQ_t)$:
\begin{equation*}
\begin{array}{ll}
\mQ_{t+1} & = \qr\Big(\left(\mQ_t + \gamma_t \mQ_t\mQ_t^\top \nabla F(\mQ_t)\right) +\gamma_t\left(\id_n -\mQ_t\mQ_t^\top \right) \nabla F(\mQ_t) \Big)\\
& = \widetilde{\mQ_t}+ \D \qr(\widetilde{\mQ_t})[\gamma_t\left(\id_n -\mQ_t\mQ_t^\top \right) \nabla F(\mQ_t)] \\
&~~ + \O\left(\gamma_t^2\|\left(\id_n -\mQ_t\mQ_t^\top \right) \nabla F(\mQ_t) \|^2\right),
\end{array}
\end{equation*}
where 
\begin{equation}
\label{eqn: qr factor}
\widetilde{\mQ_t} \mR= \mQ_t + \gamma_t \mQ_t\mQ_t^\top \nabla F(\mQ_t)\in \Real^{n \times r}
\end{equation}
denotes the QR factorization. (Note this quantity is full-rank according to \textbf{C3} in \cref{subsec: 1st convergence} and $\qr: \st(n,r) \to \st(n,r)$ is $C^\infty$ at full-rank inputs \cite{WLL12}. ) 

By the derivative formula of the QR factorization given in \cite[Algorithm 1]{Z20} (originally \cite{WLL12}), 
\begin{equation}
\label{eqn: qr taylor}
\begin{array}{ll}
\mQ_{t+1} &  = \widetilde{\mQ_t} + \gamma_t \left(\id_n - \widetilde{\mQ_t} \widetilde{\mQ_t}^\top \right)\left(\id_n - \mQ_t\mQ_t^\top \right)\nabla F(\mQ_t) \mR^{-1}\\
& ~~+\widetilde{\mQ_t}\,s\left( \gamma_t\widetilde{\mQ_t}^\top  \left(\id_n - \mQ_t\mQ_t^\top \right)\nabla F(\mQ_t) \mR^{-1}\right)\\
&~~+ \O\left(\gamma_t^2\|\left(\id_n -\mQ_t\mQ_t^\top \right) \nabla F(\mQ_t) \|^2\right)\\
&=\widetilde{\mQ_t}+ \gamma_t \left(\id_n - \mQ_t\mQ_t^\top \right)\nabla F(\mQ_t)\mR^{-1}+\O\left(\gamma_t^2\|\left(\id_n -\mQ_t\mQ_t^\top \right) \nabla F(\mQ_t) \|^2\right).
\end{array}
\end{equation}
To see the second equality, note  $\id_n - \widetilde{\mQ_t} \widetilde{\mQ_t}^\top \! = \id_n - \mQ_t\mQ_t^\top$, because $\widetilde{\mQ_t} $ and $\mQ_t$ are both the orthonormal bases for $\colspan(\mQ_t).$ Moreover, given the linear operator $s : \Real^{r \times r} \to \Real^{r \times r}$,\footnote{The operator $s$ is a composition of element-wise multiplication and skew-symmetrization. See \cite[Algorithm 1]{Z20} for details.} the term $s\left( \gamma_t\widetilde{\mQ_t}^\top  \left(\id_n - \mQ_t\mQ_t^\top \right)\nabla F(\mQ_t) \mR^{-1}\right) = \bf 0$ as $\widetilde{\mQ_t}^\top  \left(\id_n - \mQ_t\mQ_t^\top \right) = \bf 0.$

From \eqref{eqn: qr taylor}, 
\begin{equation}
\label{eqn: P taylor}
\begin{array}{ll}
\mP_{t+1} - \mP_t & = \mQ_{t+1}\mQ_{t+1}^\top - \widetilde{\mQ_t} \widetilde{\mQ_t}^\top\!\!\!\!\! \\
& =   2\gamma_t \sym\Big(\left(\id_n - \mQ_t\mQ_t^\top \right)\nabla F(\mQ_t) \mR^{-1} \widetilde{\mQ_t}^\top\Big)  \\
&~~+\O\left(\gamma_t^2\left\|\left(\id_n -\mQ_t\mQ_t^\top \right) \nabla F(\mQ_t)\mR^{-1} \right\|^2\right)\\
& ~~+  \O\left(\gamma_t^2\left\|\sym\Big(\left(\id_n -\mQ_t\mQ_t^\top \right) \nabla F(\mQ_t)\mQ_t^\top\Big) \right\|^2\right).
\end{array}
\end{equation}
Here the second line refers to the symmetric term \[\gamma_t^2\left(\id_n -\mQ_t\mQ_t^\top \right) \nabla F(\mQ_t)\mR^{-1}\mR^{-\top}\nabla F(\mQ_t)^\top\left(\id_n -\mQ_t\mQ_t^\top \right),\] whose norm is less than $\gamma_t^2\left\|\left(\id_n -\mQ_t\mQ_t^\top \right) \nabla F(\mQ_t)\mR^{-1} \right\|^2$.

We transform the right-hand side of \eqref{eqn: P taylor} from an expression in terms of $\mQ$ to one with $\mP$.

For the first component of \eqref{eqn: P taylor}, we separate $\left(\id_n - \mQ_t\mQ_t^\top \right)\nabla F(\mQ_t) \mR^{-1} \widetilde{\mQ_t}^\top$ into two parts:
\[
\texttt{Factor~1} := \left(\id_n - \mQ_t\mQ_t^\top \right)\nabla F(\mQ_t){\color{blue} \mQ_t^\top}, ~~~~\texttt{Factor~2}  := {\color{blue} \mQ_t} \mR^{-1} \widetilde{\mQ_t}^\top,
\]
using $\mQ_t^\top \mQ_t = \id_r$. By \eqref{eqn: euclid grad FQ}, 
\begin{equation}
\label{eqn: B: factor 1}
\texttt{Factor~1} = 2(\id_n - \mP_t) w(\mP_t) \mQ_t \mQ_t^\top =  2(\id_n - \mP_t)w(\mP_t) \mP_t.
\end{equation}
Also,
\begin{equation}
\label{eqn: B: factor 2}
\texttt{Factor~2}  = \mQ_t \mR^{-1} \widetilde{\mQ_t}^\top = (\mQ_t^\top)^\dag \mR^{-1} (\widetilde{\mQ_t})^\dag = \left( \widetilde{\mQ_t}\mR\mQ_t^\top\right)^\dag.
\end{equation}
by the fact that pseudoinversion satisfies $\mX^\dag\mY^\dag = (\mY\mX)^\dag$ if $\mX$ and $\mY$ respectively have orthonormal rows and columns. Inserting \eqref{eqn: qr factor} gives
\[
\texttt{Factor~2} = \left( \mQ_t\mQ_t^\top + \gamma_t \mQ_t\mQ_t^\top \nabla F(\mQ_t) \mQ_t^\top \right)^\dag = \left( \mP_t + 2\gamma_t \mP_t w(\mP_t) \mP_t\right)^\dag =  \alpha_t(\mP_t)^\dag,
\]
where the last equality is by \eqref{eqn: euclid grad FQ} and the definition of $\alpha_t$ \eqref{eqn: alpha}. 
Moreover, 
\[
\texttt{Factor~1}^\top \times \texttt{Factor~2} = \mQ_t \nabla F(\mQ_t)^\top \underbrace{\left(\id_n - \mQ_t\mQ_t^\top \right) \mQ_t }_{\bf 0}\mR^{-1} \widetilde{\mQ_t}^\top = \bf 0.
\]
Combining \eqref{eqn: B: factor 1} and \eqref{eqn: B: factor 2}, $\left(\id_n - \mQ_t\mQ_t^\top \right)\nabla F(\mQ_t) \mR^{-1} \widetilde{\mQ_t}^\top$ equals:
\begin{equation}
\label{eqn: factor 12}
\begin{array}{l}
\texttt{Factor~1} \times \texttt{Factor~2} 
 = 2(\id_n - \mP_t)w(\mP_t)\mP_t\alpha_t(\mP_t)^\dag\\ 
= (\texttt{Factor~1}+ \texttt{Factor~1}^\top) \times \texttt{Factor~2} =
2\sym\Big(2(\id_n - \mP_t)w(\mP_t) \mP_t \Big)\alpha_t(\mP_t)^\dag\\
=2\grad f(\mP_t)\alpha_t(\mP_t)^\dag.
\end{array}
\end{equation}
The formula of $\grad f$ \eqref{eqn: riem grad fP} is used in the last equality. Therefore, the first component of \eqref{eqn: P taylor} reads
\begin{equation}
\label{eqn: first line}
2\gamma_t \sym\Big(\left(\id_n - \mQ_t\mQ_t^\top \right)\nabla F(\mQ_t) \mR^{-1} \widetilde{\mQ_t}^\top\Big) = 4\gamma_t\sym\Big(\grad f(\mP_t)\alpha_t(\mP_t)^\dag \Big).
\end{equation}

As for second component of \eqref{eqn: P taylor}, note
\begin{equation*}
\begin{array}{ll}
\left\|\left(\id_n -\mQ_t\mQ_t^\top \right) \nabla F(\mQ_t)\mR^{-1}\right\|^2&=\left\|\left(\id_n -\mQ_t\mQ_t^\top \right) \nabla F(\mQ_t)\mR^{-1}\widetilde{\mQ_t}^\top \right\|^2\\
&=  4\left\|\grad f(\mP_t)\alpha_t(\mP_t)^\dag \right\|^2.
\end{array}
\end{equation*}
The orthogonality of $\widetilde{\mQ_t}$, which preserves norms, gives the first equality. The second is from \eqref{eqn: factor 12}. 

The third component of \eqref{eqn: P taylor} is dealt with directly by \eqref{eqn: grad PQ},
\begin{equation*}
\begin{array}{ll}
\left\|\sym\Big(\left(\id_n -\mQ_t\mQ_t^\top \right) \nabla F(\mQ_t)\mQ_t^\top\Big) \right\|^2&
 = \left\| \grad f(\mP_t) \right\|^2.
\end{array}
\end{equation*}

Thus, \eqref{eqn: P taylor} consists of the main component \eqref{eqn: first line} and the residuals 

$\O \left( \gamma_t^2 \left\|\grad f(\mP_t)\alpha_t(\mP_t)^\dag\right\|^2 \right)$, $\O \left( \gamma_t^2\left\|\grad f(\mP_t)\right\|^2 \right)$ as needed in \eqref{eqn: difference}. This completes the proof of \cref{difference}.
\end{proof}

\section{Proof of \cref{alpha}}
\label{proof alpha}
\begin{proof}
Consider the step size $\gamma_t > 0$ and  projector $\mP\in \gr(n,r).$ 
Definition \eqref{eqn: alpha} reads
\begin{equation}
\label{eqn: alpha P}
\alpha_t(\mP) = \mP + 2\gamma_t \mP w(\mP) \mP = \mP\left(\id_n + 2\gamma_t w(\mP) \right) \mP.
\end{equation}
So $\colspan(\alpha_t(\mP)) \subseteq \colspan(\mP)$. Also for any vector $\vx \in \colspan(\mP)$, 
\[
\vx^\top \alpha_t(\mP) \vx = \vx^\top \mP (\id_n + 2\gamma_t w(\mP) )  \mP\vx = \vx^\top(\id_n + 2\gamma_t w(\mP) )\vx \geq \|\vx\|^2,
\]
as $ \mP \vx = \vx$ and $\id_n + 2\gamma_t w(\mP) \succeq \id_n.$
It follows
$\colspan(\alpha_t(\mP)) = \colspan(\mP)$, and
$\alpha_t(\mP)$ has only $r$ non-zero singular values, each   greater than $1$. 
By properties of pseudoinverses,
\begin{equation}
\label{eqn: alpha inv}
\alpha_t(\mP)^\dag\alpha_t(\mP) = \textbf{Proj}_{\text{range}(\alpha_t(\mP))} = \mP,
\end{equation}
and $\alpha_t(\mP)^{\dagger}$ has only $r$ non-zero singular values, each less than 1. Hence
\begin{equation}
\label{eqn: alpha inv norm}
\|\alpha_t(\mP)^\dag\| \leq \sqrt{r}.
\end{equation}

Constrain the step size by
\begin{equation}
\label{eqn: gamma 1}
0 < \gamma_t \leq \Gamma_1 := \frac{1}{4\sqrt{r}\sup_{\mP \in\gr(n,r)}\|\mP w(\mP) \mP\|},
\end{equation}
where $\sup_{\mP \in\gr(n,r)}\|\mP w(\mP) \mP\|$ is finite because the function $\mP \mapsto \mP w(\mP) \mP$ is continuous  and $\gr(n,r)$ is compact. 
(We  exclude the trivial case  $ \mP w(\mP)\mP \equiv 0$ when $\tX = 0$.)
From \eqref{eqn: alpha inv}, 
\begin{equation}
\label{eqn: alpha trick}
\mP - \alpha_t(\mP)^\dag =  \alpha_t(\mP)^\dag \alpha_t(\mP)- \alpha_t(\mP)^\dag\mP = - \alpha_t(\mP)^\dag(\mP -\alpha_t(\mP)).
\end{equation}
Then by \eqref{eqn: alpha inv norm} and \eqref{eqn: gamma 1},
\begin{equation}
\label{eqn: P - alpha norm}
\begin{array}{ll}
\| \mP - \alpha_t(\mP)^\dag\| &= \| -\alpha_t(\mP)^\dag\left(\mP-\alpha_t(\mP) \right)\|\\
& \leq \|\alpha_t(\mP)^\dag\| \|\mP-\alpha_t(\mP)\| \\
& \displaystyle \leq 2\gamma_t\sqrt{r} \|\mP w(\mP) \mP\| \leq  \frac{1}{2}. 
\end{array}
\end{equation}

We now compare $\| \grad f(\mP)\|$ and $\left\|\grad f(\mP)\alpha_t(\mP)^\dag\right\|$. By \eqref{eqn: grad PQ}, 
\[
\grad f(\mP)\mP = (\id_n -\mP) w(\mP) \mP.
\] 
Inserting the above equation  into \eqref{eqn: grad PQ}, 
\[
\grad f(\mP)=\sym\Big(2(\id_n -\mP) w(\mP) \mP\Big)=\sym\Big(2\grad f(\mP)\mP\Big).
\]
Since 
\[
\left\langle \grad f(\mP)\mP, \mP\grad f(\mP)\right\rangle = \tr\left( (\id_n -\mP) w(\mP) \underbrace{\mP(\id_n -\mP)}_{\bf 0} w(\mP) \mP\right)=0,
\] 
we have 
\begin{equation}
\label{eqn: grad fP P}
\left\| \grad f(\mP)\right\| =\sqrt{2} \left\|\grad f(\mP)\mP\right\|.
\end{equation}
Then due to triangle inequality and \eqref{eqn: alpha trick}, 
\begin{equation*}
\begin{array}{ll}
\| \grad f(\mP)\| & \leq  \sqrt{2}  \left\|\grad f(\mP)\alpha_t(\mP)^\dag\right\| +  \sqrt{2}  \left\|\grad f(\mP)\left(\mP - \alpha_t(\mP)^\dag\right)\right\|\\
& \displaystyle = \sqrt{2} \left\|\grad f(\mP)\alpha_t(\mP)^\dag\right\| +  \sqrt{2}  \left\|\grad f(\mP)\alpha_t(\mP)^\dag\left(\mP - \alpha_t(\mP)^\dag\right)\right\|\\
& \displaystyle\leq \sqrt{2} \left\|\grad f(\mP)\alpha_t(\mP)^\dag\right\| \left(1+\left\|\mP - \alpha_t(\mP)^\dag  \right\|\right) \\
& \displaystyle= \frac{3\sqrt{2}}{2}\left\|\grad f(\mP)\alpha_t(\mP)^\dag\right\|.
\end{array}
\end{equation*}
Inserting \eqref{eqn: P - alpha norm} shows the last inequality. Thus,
\[
\left\|\grad f(\mP)\alpha_t(\mP)^\dag\right\| \geq \frac{\sqrt{2}}{3} \left\|\grad f(\mP)\right\|.
\]

On the other hand, as $\alpha_t(\mP)^\dag = \mP\alpha_t(\mP)^\dag,$ it holds
\begin{equation*}
\begin{array}{ll}
\left\|\grad f(\mP)\alpha_t(\mP)^\dag\right\|
& \leq  \left\|\grad f(\mP)\mP\right\|  \left\|\alpha_t(\mP)^\dag\right\|\\
& \leq  \displaystyle\frac{\sqrt{2r}}{2} \left\|\grad f(\mP)\right\|.
\end{array}
\end{equation*}
The results \eqref{eqn: grad fP P} and \eqref{eqn: alpha inv norm} give the last derivation. Combining the above two inequalities completes the proof of \cref{alpha}. 
\end{proof}

\section{Proof of \cref{hess Dphi formula} }
\label{proof hess Dphi formula}
\begin{proof}
Let $\mP\in\gr(n,r)$ and $\Delta \mP \in \tangent_{\mP}\gr(n,r)$. We first show the computation of the Riemannian Hessian. Applying \cite[Equation (2.109)]{HHT07}, 
\begin{equation}
\label{eqn: riem hess f 2}
\begin{array}{l}
\hess f(\mP)[\Delta \mP] \\
= \ad_\mP^2\left(\nabla^2f(\mP)[\Delta \mP] \right)-\ad_\mP\ad_{\nabla f(\mP)}(\Delta \mP)\\
= \ad_\mP\Big(\mP \nabla^2 f(\mP)[\Delta \mP]- (\nabla^2 f(\mP)[\Delta \mP])\mP - \nabla f(\mP) \Delta \mP + \Delta \mP\nabla f(\mP) \Big)\\
= \mP(\nabla^2 f(\mP)[\Delta \mP])(\id_n -\mP) + (\id_n - \mP)(\nabla^2 f(\mP) [\Delta \mP])\mP\\
 ~~-\mP\nabla f(\mP)\Delta \mP + \mP \Delta \mP \nabla f(\mP) + \nabla f(\mP)(\Delta \mP) \mP-\Delta \mP \nabla f(\mP)\mP\\
 = 2\sym\Big((\id_n - \mP)(\nabla^2 f(\mP) [\Delta \mP])\mP +  \nabla f(\mP) (\Delta \mP) \mP-\Delta \mP  \nabla f(\mP)\mP\Big). 
\end{array}
\end{equation}

Recall \eqref{eqn: euclid grad f} that $\nabla f(\mP) = w(\mP)\mP+\mP w(\mP)$.
We differentiate $\nabla f(\mP)$ to obtain the Euclidean Hessian of $f$:
\begin{equation}
\label{eqn: euclid hess f}
\nabla^2f(\mP)[\Delta \mP] = \left(\D w(\mP)[\Delta \mP]\right) \mP+ w(\mP)\Delta \mP+\Delta \mP w(\mP) + \mP \left(\D w(\mP)[\Delta \mP]\right),
\end{equation}
where the differential is
\begin{equation}
\label{eqn: D w}
\begin{array}{l}
\D w(\mP) [\Delta \mP] \\
= 2d(d-1) \sym\Big(\left\langle \tX\cdot\left(\id_n, \Delta \mP, \mP, \dots, \mP \right),\tX\cdot\left(\id_n, \mP, \mP, \dots, \mP \right) \right\rangle_{-1}\Big)\\
 = d(d-1) \left\langle \tX\cdot\left(\id_n, (\Delta \mP)\mP+\mP\Delta \mP, \mP, \dots, \mP \right),\tX\cdot\left(\id_n, \id_n, \mP, \dots, \mP \right) \right\rangle_{-1}\\
 = d(d-1)\left\langle \tX\cdot\left(\id_n, \Delta \mP, \mP, \dots, \mP \right),\tX\cdot\left(\id_n, \id_n, \mP, \dots, \mP \right) \right\rangle_{-1}.
\end{array}
\end{equation}
The last equality follows from  $(\Delta \mP)\mP+\mP\Delta \mP = \Delta \mP$ in \eqref{eqn: tangent space}.

For the Riemannian Hessian \eqref{eqn: riem hess f 2}, we first consider the term involving 

$\nabla^2 f(\mP) [\Delta \mP]$: 
\begin{equation}
\label{eqn: riem hess term 1}
\begin{array}{l}
 (\id_n - \mP)(\nabla^2 f(\mP) [\Delta \mP])\mP\\
=(\id_n - \mP)\left( \D w(\mP) [\Delta \mP]\right)\mP+(\id_n - \mP)w(\mP) \left(\Delta \mP \right) \mP+(\id_n - \mP)\Delta \mP w(\mP)\mP\\
= (\id_n - \mP)v(\mP) [\Delta \mP]+ (\id_n - \mP)w(\mP) \left(\Delta \mP \right) \mP+(\id_n - \mP)\Delta \mP w(\mP)\mP.
\end{array}
\end{equation}
The last equality is due to $v(\mP)[\Delta \mP] = \left(\D w(\mP)[\Delta \mP]\right) \mP$ by \eqref{eqn: v} and \eqref{eqn: D w}.

We then compute the remaining terms in \eqref{eqn: riem hess f 2}, which involve $\nabla f(\mP)$:  
\begin{equation}
\label{eqn: riem hess term 2}
\begin{array}{l}
\nabla f(\mP)(\Delta \mP) \mP-\Delta \mP  \nabla f(\mP)\mP\\
 = \left(w(\mP)\mP+\mP w(\mP) \right)\left(\Delta \mP \right) \mP-\Delta \mP\left(w(\mP)\mP+\mP w(\mP) \right)\mP\\
=\mP w(\mP)\left(\Delta \mP \right)\mP -\Delta \mP w(\mP) \mP - \left(\Delta \mP \right)\mP w(\mP)\mP.
\end{array}
\end{equation}
In the last equality, we eliminate the term that includes $\mP( \Delta \mP) \mP$, since this is $\bf 0$ due to \eqref{eqn: tangent space}. 

We add  \eqref{eqn: riem hess term 1} and \eqref{eqn: riem hess term 2} to obtain
\begin{equation*}
\begin{array}{l}
\hess f(\mP)[\Delta \mP] \\
=2\sym\Big((\id_n - \mP) v(\mP)[\Delta \mP]  + (\id_n - \mP+\mP)w(\mP) \left(\Delta \mP \right)\mP\\ 
 +  (\id_n - \mP - \id_n) \left( \Delta \mP\right) w(\mP)\mP- \left(\Delta \mP \right)\mP w(\mP)\mP\Big)\\
=2\sym\Big( (\id_n - \mP) v(\mP)[\Delta \mP] +w(\mP) \left(\Delta \mP \right)\mP - (\mP \Delta \mP+  \left( \Delta \mP\right)\mP) w(\mP)\mP\Big)\\
 = 2\sym\Big((\id_n - \mP) v(\mP)[\Delta \mP] +w(\mP) \left(\Delta \mP \right)\mP- \Delta \mP w(\mP)\mP\Big).
\end{array}
\end{equation*}
This is the formula for the Riemannian Hessian in \eqref{eqn: riem hess f}. 

Next, we compute the differential $\D\phi(\mP)[\Delta \mP].$ This is the term that is linear in $\Delta \mP$ in the Taylor expansion of the difference: 
\begin{equation}
\label{eqn: D phi diff}
 \phi(\mP+ \Delta \mP)- \phi(\mP) = \Phi(\mQ+\Delta \mQ)\Phi(\mQ+\Delta \mQ)^\top- \Phi(\mQ) \Phi(\mQ)^\top.
\end{equation}

So we start with $\Phi(\mQ+\Delta\mQ)$ and $\Phi(\mQ)$. From the first-order Taylor expansion of $\qr$, 
\begin{equation*}
\begin{array}{ll}
\Phi(\mQ+\Delta \mQ) &= \qr\left(\mQ+\Delta\mQ+\gamma\nabla F(\mQ+\Delta\mQ)\right)\\
& = \qr\Big(\left(\mQ+\gamma\nabla F(\mQ) \right) + \Delta\mQ+\gamma\left(\nabla F(\mQ+\Delta\mQ) - \nabla F(\mQ)\right)\Big)\\
& = \qr\Big(\left(\mQ+\gamma\nabla F(\mQ) \right) + \Delta\mQ+\gamma\nabla^2F(\mQ)[\Delta \mQ]+\O\left(\|\Delta \mQ\|^2\right)\Big)\\
& = \qr\left(\mQ+\gamma\nabla F(\mQ) \right) \\
&~~+ \D\qr\left(\mQ+\gamma\nabla F(\mQ) \right)\left[\Delta\mQ+\gamma\nabla^2 F(\mQ)[\Delta \mQ]+\O\left(\|\Delta \mQ\|^2\right)\right].
\end{array}
\end{equation*}
Write $\Phi(\mQ) \mR$ for the QR decomposition of $\mQ+\gamma \nabla F(\mQ)$. We again apply the differentiation formula for $\qr$ in \cite[Algorithm 1]{Z20} and get
\begin{equation}
\label{eqn: Phi(Q+delta)}
\begin{array}{ll}
\Phi(\mQ+\Delta \mQ) 
& = \Phi(\mQ) + \left(\id_n - \Phi(\mQ)\Phi(\mQ)^\top \right)\left(\Delta \mQ+\gamma\nabla^2 F(\mQ)[\Delta \mQ]\right)  \mR^{-1}\\
&+ \Phi(\mQ) \mS+ \O\left(\|\Delta \mQ\|^2\right).
\end{array}
\end{equation}
Here $\mS \in \Real^{r \times r}$ is a  skew-symmetric matrix that is linear in $\Delta \mQ$. 
The Euclidean gradient $\nabla F(\mQ)$ is given by \eqref{eqn: euclid grad FQ}, and the Euclidean Hessian of $F$ is calculated as 
\begin{equation}
\label{eqn: euclid hess FQ}
\begin{array}{l}
\nabla^2 F(\mQ)[\Delta \mQ] \\
=2d\left\langle \tX \cdot (\id_n, \mQ, \dots, \mQ), \tX \cdot (\id_n, \mQ, \dots, \mQ) \right\rangle_{-1}\Delta \mQ\\
~~+2d(d-1)\left\langle \tX \cdot (\id_n,(\Delta \mQ) \mQ + \mQ\Delta \mQ, \mQ, \dots, \mQ), \tX \cdot (\id_n, \id_n, \mQ, \dots, \mQ) \right\rangle_{-1}\mQ\\
 = 2w(\mP)\Delta \mQ+2d(d-1)\left\langle \tX \cdot (\id_n,\Delta\mP, \mQ, \dots, \mQ), \tX \cdot (\id_n, \id_n, \mQ, \dots, \mQ) \right\rangle_{-1}\mQ,
\end{array}
\end{equation}
due to \eqref{eqn: w} and  $\Delta \mP = \left(\Delta \mQ\right) \mQ^\top + \mQ \left(\Delta \mQ\right)^\top.$

We substitute \eqref{eqn: Phi(Q+delta)} into \eqref{eqn: D phi diff}. We only keep the terms that are linear in $\Delta \mQ$ and get
\begin{equation}
\label{eqn: D phi expand}
\begin{array}{ll}
\D \phi(\mP) [\Delta \mP]  & = 2\,\sym\Big(\left(\id_n - \Phi(\mQ)\Phi(\mQ)^\top \right)\left(\Delta \mQ+\gamma\nabla^2 F(\mQ)[\Delta \mQ] \right)\mR^{-1}\Phi(\mQ)^\top\Big) \\
& =2\,\sym\Big(\left(\id_n - \phi(\mP) \right)\left(\Delta \mQ+\gamma\nabla^2 F(\mQ)[\Delta \mQ] \right)\left(\mR^{-1}\Phi(\mQ)^\top\right)\Big)\\
& = 2\,\sym\Big(\left(\id_n - \phi(\mP) \right)\left(\Delta \mQ+\gamma\nabla^2 F(\mQ)[\Delta \mQ] \right)\left(\mQ+\gamma \nabla F(\mQ)\right)^\dag\Big).
\end{array}
\end{equation}
For the first equality, we use  $\Phi(\mQ)\mS\Phi(\mQ)^\top + \Phi(\mQ)\mS^\top\Phi(\mQ)^\top = \bf 0$ and the fact  $\mS$ is linear in $\Delta \mQ$ 

It remains to transform the differential expression \eqref{eqn: D phi expand} in $\mQ, \Delta \mQ$ to be in terms of $\mP, \Delta \mP$. To do this, we split \eqref{eqn: D phi expand} into two factors:
\begin{equation*}
\begin{array}{l}
\texttt{Factor~1}:= \left(\id_n - \phi(\mP) \right)\left(\Delta \mQ+\gamma\nabla^2 F(\mQ)[\Delta \mQ] \right){\color{blue} \mQ^\top}, \\
\texttt{Factor~2} :={\color{blue} \mQ}\left(\mQ+\gamma \nabla F(\mQ)\right)^\dag. 
\end{array}
\end{equation*}

We first tackle \texttt{Factor~2}. As $\mQ^\top$ and $\mQ+\gamma \nabla F(\mQ)$ respectively have orthonormal rows and columns, we have 
\begin{equation*}
\mQ\left(\mQ+\gamma \nabla F(\mQ)\right)^\dag = \left(\mQ^\top\right)^\dag\left(\mQ+\gamma \nabla F(\mQ)\right)^\dag = \left( (\mQ+ \gamma\nabla F(\mQ))\mQ^\top\right)^\dag.
\end{equation*}
Inserting \eqref{eqn: euclid grad FQ} into the above equation, 
\begin{equation}
\label{eqn: factor 2}
\begin{array}{ll}
 \texttt{Factor~2} & =\left(\left(\mQ+ 2\gamma w(\mP)\mQ\right) \mQ^\top\right)^{\dag}= \left(\mP +2\gamma w(\mP)\mP\right)^{\dag}.
\end{array}
\end{equation}

Next, we consider \texttt{Factor~1}. By \eqref{eqn: euclid hess FQ}, it equals 
\begin{equation}
\label{eqn: factor 1 unfinished}
\begin{array}{l}
\left(\id_n - \phi(\mP) \right)\left(\Delta \mQ+\gamma\nabla^2 F(\mQ)[\Delta \mQ]\right) \mQ^\top\\
= \left(\id_n -  \phi(\mP)\right) \left(\id_n + 2\gamma w(\mP)\right)(\Delta \mQ)\mQ^\top  \\
+\left(\id_n -  \phi(\mP) \right) \left(2\gamma d(d-1)\left\langle \tX \cdot (\id_n,\Delta\mP , \mQ, \dots, \mQ), \tX \cdot (\id_n, \id_n, \mQ, \dots, \mQ) \right\rangle_{-1}\mQ\mQ^\top  \right)\\
= \left(\id_n -  \phi(\mP) \right)\left(\left(\id_n + 2\gamma w(\mP) \right)(\Delta \mP - \mQ(\Delta \mQ)^\top) +2 \gamma v(\mP)[\Delta \mP]\right).
\end{array}
\end{equation}
The last equality is from the definition of $v(\mP)[\Delta \mP]$ \eqref{eqn: v}.

However, note that
\begin{equation*}
\begin{array}{l}
 \left(\id_n - \phi(\mP) \right)\left(\id_n + 2\gamma w(\mP) \right)\mQ(\Delta \mQ)^\top  \\
 =  \left(\id_n -  \phi(\mP) \right)\left( \mQ + 2\gamma w(\mP)\mQ\right)(\Delta\mQ)^\top\\
  =  \left(\id_n -  \Phi(\mQ)\Phi(\mQ)^\top \right)\left( \mQ + \gamma \nabla F(\mQ)\right)(\Delta\mQ)^\top = \bf 0.
  \end{array}
\end{equation*}
Here, we replace $2w(\mP)\mQ$ by $\nabla F(\mQ)$ by \eqref{eqn: euclid grad FQ}.  
Then from $\colspan(\mQ + \gamma \nabla F(\mQ))$ equals $\colspan(\Phi(\mQ))$ by \eqref{eqn: Q update}, it holds $\left(\id_n -  \Phi(\mQ)\Phi(\mQ)^\top \right)\left( \mQ + \gamma \nabla F(\mQ)\right)=\bf 0$.

Returning to \eqref{eqn: factor 1 unfinished}, 
\begin{equation}
\label{eqn: factor 1}
\begin{array}{ll}
\texttt{Factor~1} &= \left(\id_n -  \phi(\mP) \right)\left(\left(\id_n +2 \gamma w(\mP) \right)(\Delta \mP) +2 \gamma v(\mP)[\Delta \mP]\right)\\
& = \left(\id_n -  \phi(\mP) \right)\left(\Delta \mP +2 \gamma w(\mP)(\Delta \mP) +2 \gamma v(\mP)[\Delta \mP]\right).
\end{array}
\end{equation}

Putting \eqref{eqn: D phi expand}, \eqref{eqn: factor 2} and \eqref{eqn: factor 1} together, we obtain the  formula for the differential $\D\phi(\mP)[\Delta \mP]$ \eqref{eqn: D phi}. The proof of \cref{hess Dphi formula} is complete. 
\end{proof}

\section{Proof of \cref{hess Dphi}}
\label{proof hess Dphi}
\begin{proof}
Let $\mP \in \gr(n,r)$ be a first-order critical point of $f$ \eqref{eqn: cost 2}. 
The Riemannian Hessian $\hess f(\mP) : \tangent_\mP \gr(n,r) \rightarrow \tangent_\mP \gr(n,r)$ is diagonalizable because it is self-adjoint.
The dimension of the tangent space $\tangent_\mP \gr(n,r) $ is $\text{dim}(\gr(n,r))=r(n-r)$.
Let $\Delta \mP_1, \dots,  \Delta \mP_{r(n-r)}\in \tangent_\mP \gr(n,r) $ be a basis of eigenmatrices of $\hess f(\mP)$ with corresponding eigenvalues $\lambda_1, \dots, \lambda_{r(n-r)},$ i.e. %for $i \in [r(n-r)]$,
\[
\hess f(\mP)[\Delta \mP_i] =\lambda_i \Delta \mP_i
\]
for $i \in [r(n-r)]$.
Assume there is at least one strictly positive eigenvalue of $\hess f(\mP)$, w.l.o.g. $\lambda_1 > 0$. 

We  analyze the eigenvalues of $\D\phi(\mP)$ by expressing the differential with respect to $\{\Delta \mP_i \}$. 
Under this basis, the differential is represented by the matrix $\mM \in \Real^{r(n-r) \times r(n-r)}$ with entries $m_{i,j}$ equaling to 
\begin{equation}
\label{eqn: m}
\left\langle 2\,\sym\Big( (\id_n - \mP) \left( \Delta \mP_i + 2\gamma w(\mP) \Delta \mP_i+ \gamma v(\mP) [\Delta \mP_i]\right) (\mP+2\gamma w(\mP)\mP)^\dag\Big), ~~\Delta \mP_j\right\rangle,
\end{equation}
for $i, j \in [r(n-r)]$. Here we apply the formula \eqref{eqn: D phi} and the fact that $\mP = \phi(\mP)$ as critical points of $f$ are fixed points of $\phi$ by \textbf{C3} in \cref{subsec: 1st convergence}. The eigenvalues of the operator $ \D \phi(\mP)$ are the same as the eigenvalues of the matrix $\mM.$

We next compute the matrix $\mM$. Due to the rank constraint on $\mP$, we utilize its block structure to facilitate the computation. 
W.l.o.g., we suppose the rank-$r$ orthogonal projection $\mP$ is:
\begin{equation*}
\mP  = 
\left[
\begin{array}{cc}
\id_{r}&{\bf 0}_{r \times (n-r)}\\
{\bf 0}_{(n-r) \times r} &{\bf 0}_{(n-r) \times (n-r)}
\end{array}
\right] \in \gr (n,r),
\end{equation*}
by choosing an appropriate for $\Real^n$.

From \eqref{eqn: riem grad fP},
$
\grad f(\mP) = \sym\Big(2 (\id_n-\mP) w(\mP) \mP\Big) = {\bf 0}, 
$
and \eqref{eqn: tangent space}, $
\Delta \mP = \left(\Delta \mP\right)\mP + \mP\Delta \mP$ for all $\Delta \mP \in  \tangent_\mP \gr(n,r)$,
we obtain block matrices for $w(\mP)$ and $\left\{\Delta \mP_i\right\}$. 
Namely there are PSD matrices $\mA \in \S^2(\Real^r), \mC \in \S^2(\Real^{n-r})$ and matrices $\mH_i \in \Real^{(n-r) \times r}$ for $i \in [r(n-r)]$ such that
\begin{equation}
\label{eqn: block}
w(\mP)  =
\left[
\begin{array}{cc}
\mA&\bf 0\\
\bf 0 &\mC
\end{array}
\right] \in \S^2(\Real^n),
~~~~
\Delta \mP_i =
\left[
\begin{array}{cc}
\bf 0&\mH_i^\top\\
\mH_i &\bf 0
\end{array}
\right] \in\tangent_{\mP}\gr (n,r).
\end{equation}

Since $\Delta \mP_i$ are eigenmatrices of $\hess f(\mP)$, by \eqref{eqn: riem hess f} we see
\begin{equation*}
\begin{array}{l}
2\,\sym\Big((\id_n - \mP) v(\mP) [\Delta \mP_i] \Big)\\
= \hess f(\mP)[\left(\Delta \mP\right)_i]   -2\, \sym\Big(w(\mP) (\Delta \mP_i) \mP-\Delta \mP_i w(\mP)\mP\Big) \\
%= 
%\left[
%\begin{array}{cc}
%\bf 0&\lambda_i\mH_i^\top\\
%\lambda_i\mH_i &\bf 0
%\end{array}
%\right]
%\\
%-2\,\sym\left(
%\left[
%\begin{array}{cc}
%\mA&\bf 0\\
%\bf 0 &\mC
%\end{array}
%\right]
%\left[
%\begin{array}{cc}
%\bf 0&\mH_i\\
%\mH_i &\bf 0
%\end{array}
%\right]
%\left[
%\begin{array}{cc}
%\id_r&\bf 0\\
%\bf 0 &\bf 0
%\end{array}
%\right]
% -
%\left[
%\begin{array}{cc}
%\bf 0&\mH_i^\top\\
%\mH_i&\bf 0
%\end{array}
%\right]
%\left[
%\begin{array}{cc}
%\mA&\bf 0\\
%\bf 0 &\mC
%\end{array}
%\right]
%\left[
%\begin{array}{cc}
%\id_r&\bf 0\\
%\bf 0 &\bf 0
%\end{array}
%\right]\right)\\\\
 = 
\left[
\begin{array}{cc}
\bf 0&\lambda_i\mH_i^\top -\mH_i^\top\mC+ \mA\mH_i^\top \\
\lambda_i\mH_i-\mC\mH_i+\mH_i\mA& \bf 0
\end{array}
\right].
 \end{array}
\end{equation*}
It follows that 
\begin{equation*}
(\id_n -\mP) v(\mP) [\Delta \mP_i] = 
\left[
\begin{array}{cc}
\bf 0&\bf 0\\
\lambda_i\mH_i-\mC\mH_i+\mH_i\mA& \bf 0
\end{array}
\right],
\end{equation*}
because $(\id_n - \mP)v(\mP) [\Delta \mP_i]$ \eqref{eqn: v} is in the format $(\id_n - \mP) (\cdot) \mP$ so only its left bottom block is nonzero.

To get $m_{i,j}$ \eqref{eqn: m}, we compute
\begin{equation*}
\begin{array}{l}
(\id_n - \mP)\left(\Delta \mP_i  + 2\gamma w(\mP)\Delta \mP_i +2\gamma v(\mP)[\Delta \mP_i]  \right) \\ 
= 
\left[
\begin{array}{cc}
\bf 0&\bf 0\\
\bf 0& \id_{n-r}
\end{array}
\right]
\left[
\begin{array}{cc}
\id_r + 2\gamma\mA&\bf 0\\
\bf 0& \id_{n-r}+ 2\gamma\mC
\end{array}
\right]
\left[
\begin{array}{cc}
\bf 0&\mH_i^\top\\
\mH_i& \bf 0
\end{array}
\right]\\
+
\left[
\begin{array}{cc}
\bf 0&\bf 0\\
2\gamma(\lambda_i\mH_i-\mC\mH_i+\mH_i\mA)& \bf 0
\end{array}
\right]\\\\
= 
\left[
\begin{array}{cc}
\bf 0&\bf 0\\
(1+2\gamma\lambda_i) \mH_i+2\gamma\mH_i\mA& \bf 0
\end{array}
\right].
\end{array}
\end{equation*}

Also,
\begin{equation*}
\left( \mP + 2\gamma w(\mP) \mP\right)^{\dag} 
 = \left(
 \left[
\begin{array}{cc}
\id_r &\bf 0\\
\bf 0& \bf 0
\end{array}
\right]
+
\left[
\begin{array}{cc}
2\gamma\mA&\bf 0\\
\bf 0&\bf 0
\end{array}
\right]
 \right)^{\dag}
 =
 \left[
\begin{array}{cc}
(\id_r+2\gamma\mA)^{-1}&\bf 0\\
\bf 0&\bf 0
\end{array}
\right].
\end{equation*}
Here we note that since $\mA$ is PSD and $\gamma > 0$, the inverse $(\id_r+2\gamma\mA)^{-1}$ exists and is symmetric. 

Thus the left component in the inner product  \eqref{eqn: m} is 
\begin{equation*}
\begin{array}{l}
2\,\sym\left( \left[
\begin{array}{cc}
\bf 0&\bf 0\\
(1+2\gamma\lambda_i) \mH_i+2\gamma\mH_i\mA& \bf 0
\end{array}
\right]
\left[
\begin{array}{cc}
(\id_r+2\gamma\mA)^{-1}&\bf 0\\
\bf 0&\bf 0
\end{array}
\right] \right)\\\\
=
\left[
\begin{array}{cc}
\bf 0&\mH_i^\top + 2\gamma\lambda_i(\id_r+2\gamma\mA)^{-1}\mH_i^\top\\
\mH_i+2\gamma\lambda_i\mH_i(\id_r+2\gamma\mA)^{-1}& \bf 0
\end{array}
\right]
\end{array}
\end{equation*}

Finally \eqref{eqn: m} reads $m_{i,j}$ equals to
\begin{equation}
\label{eqn: m value}
\begin{array}{l}
%m_{i,j} 
%& = 
\left\langle
\left[
\begin{array}{cc}
\bf 0&\mH_i^\top + 2\gamma\lambda_i(\id_r+2\gamma\mA)^{-1}\mH_i^\top\\
\mH_i+2\gamma\lambda_i\mH_i(\id_r+2\gamma\mA)^{-1}& \bf 0
\end{array}
\right],
\left[
\begin{array}{cc}
\bf 0&\mH_j^\top \\
\mH_j& \bf 0
\end{array}
\right]
 \right\rangle\\\\
  = 2 \,\left\langle \mH_i+2\gamma\lambda_i\mH_i(\id_r+2\gamma\mA)^{-1}, \mH_j \right\rangle \\\\
  = \left\{
 \begin{array}{lc}
4 \gamma\lambda_i\left \| \mH_i (\id_r+2\gamma\mA)^{-\frac{1}{2}} \right\|^2 + 1 & i = j\\
  4\gamma\lambda_i \left\langle \mH_i (\id_r+2\gamma\mA)^{-\frac{1}{2}},  \mH_j (\id_r+2\gamma\mA)^{-\frac{1}{2}}  \right\rangle & i \neq j,
 \end{array}\right.
\end{array} 
\end{equation}
where we use
\begin{equation*}
\left\langle \mH_i, \mH_j \right\rangle = \frac{1}{2}\,\left\langle \Delta\mP_i, \Delta\mP_j \right\rangle  = 
\left\{
 \begin{array}{lc}
\frac{1}{2}, & i = j\\
  0,  & i \neq j.
 \end{array}\right.
\end{equation*}

We set 
\begin{equation}
\label{eqn: L Lambda}
\begin{array}{ll}
\mL &= 
\left[
\begin{array}{c}
\text{vec}\left( \mH_1 (\id_r+2\gamma\mA)^{-\frac{1}{2}} \right)^\top\\
\vdots\\
\text{vec}\left( \mH_{r(n-r)} (\id_r+2\gamma\mA)^{-\frac{1}{2}} \right)^\top
\end{array}
\right] \in \Real^{r(n-r) \times r(n-r)},\\\\
\pmb{\Lambda} &= 
4\gamma\diag\left(\lambda_1, \dots, \lambda_{r(n-r)}\right) \in \S^2(\Real^{r(n-r)}).
\end{array}
\end{equation}
Here $\text{vec}$ denotes the vectorization operator.
Combining \eqref{eqn: m value} and \eqref{eqn: L Lambda}, we have a short-hand expression for $\mM$: 
\[
\mM = \pmb{\Lambda} \mL\mL^\top +\id_{r(n-r)}.
\]

The linear independence of $\{\mH_i\}$ implies the linear independence of $\{\mH_i(\id_r+2\gamma\mA)^{-\frac{1}{2}}\}.$ Then $\mL$ is invertible because its rows are linearly independent.
Hence $\pmb{\Lambda}\mL\mL^\top$ is similar to $\mL^\top\pmb{\Lambda}\mL = \mL^\top \left({\pmb{\Lambda}} \mL\mL^\top\right) \mL^{-\top}$. 
As
$\mL^{\top} {\pmb{\Lambda}} \mL$ is symmetric, 
$\pmb{\Lambda}\mL\mL^\top$ is diagonalizable with the same eigenvalues as $\mL^\top\pmb{\Lambda}\mL.$  Then $ \pmb{\Lambda}$ and $\mL^\top \pmb{\Lambda} \mL$ are congruent. 
So their numbers of positive eigenvalues match. 
But there is a strictly positive eigenvalue of $\pmb{\Lambda}$, namely $4\gamma\lambda_1 > 0$. 
It follows that $\mM = \pmb{\Lambda} \mL\mL^\top+\id_{r(n-r)}$ has an eigenvalue strictly greater than 1. 
The proof of \cref{hess Dphi} is complete.
\end{proof}

\section{Proof of \cref{local diff}}
\label{proof local diff}
\begin{proof}

We prove \cref{local diff} by showing that $\phi$ is $C^{\infty}$ and $\D\phi(\mP)$ is  invertible for all $\mP\in\gr(n,r)$, when $\gamma$ is small enough.

First we claim that $\Phi$ is  $C^{\infty}$ on $\st(n,r)$.
It is a composite of $\qr$ and $\mQ \mapsto \mQ + \gamma \nabla F(\mQ)$.  The latter is a polynomial map, only taking full-rank values by $\textbf{C3}$ of \cref{subsec: 1st convergence}. 
Meanwhile, $\qr$ is  $C^{\infty}$ at full-rank inputs by \cite{WLL12}. 
So $\Phi$ is $C^\infty$. 
Consequently $\phi$ is $C^\infty$ by \cite[Theorem A.27 (a)]{L12}, since $\mQ \mapsto \mQ \mQ^{\top}$ is a quotient map from $\st(n,r)$ to $\gr(n,r)$.

Next we show  $\D\Phi(\mQ)$ is  invertible for all $\mQ \in \st(n,r)$. Define  $h: \st(n,r) \times \Real_{\geq 0} \to \Real_{\geq 0}$ by
\[
h(\mQ, \gamma) := \sigma_{\min} \left(\D\Phi(\mQ)\right) = \sigma_{\min} \left(\D\qr(\mQ+\gamma\nabla F(\mQ))\right),
\]
where $\sigma_{\min}$ denotes the smallest singular value.
Here $\qr$ is $C^{\infty}$ at full-rank inputs, $\Phi(\mQ)$ is full-rank for $\gamma \geq 0$ and $\Phi(\mQ)$ is polynomial in $\mQ, \gamma$. 
It follows that
$\D\qr(\mQ+\gamma\nabla F(\mQ))$ is jointly continuous in $\mQ, \gamma$.
As the smallest singular value is also continuous, $h$ is jointly continuous in $\mQ, \gamma$. 
When $\gamma = 0,$ $\Phi$ becomes the identity map on $\st(n,r)$, so $h(\cdot, 0) \equiv 1$.  Because the Stiefel manifold is compact, there exists a constant $\Gamma^{**}(\tX,r) > 0$ such that for all $0< \gamma \leq\Gamma^{**}$ it holds $h(\cdot, \gamma) > 0$. W.l.o.g. assume $\Gamma^{**}$ is less than $\Gamma^*$ in \cref{1st convergence}.
Then $\D \Phi: \tangent_{\mQ} \st(n,r) \to \tangent_{\phi(\mQ)} \st(n,r)$ 
is invertible for all $\mQ \in \st(n,r)$ when $0< \gamma \leq \Gamma^{**}$.

To conclude, differentiate the update diagram \cref{QP update}.  By the chain rule, we get a commutative diagram of differentials between tangent spaces. 
The ``north-then-west" path is a composition of surjections because $\D\Phi(\mQ_t)$ is surjective and $\mQ \mapsto \mQ \mQ^{\top}$ is a quotient map.
By commutativity, the ``west-then-north" composition is surjective.  So the second map there must be a surjection, i.e. $\D\phi(\mP_t) : \tangent_{\mP_t} \gr(n,r) \rightarrow \tangent_{\mP_{t+1}}\gr(n,r)$ is surjective.
Therefore it is invertible, as the tangent spaces share the same dimension.
The proof of \cref{local diff} is complete.
\end{proof}

\end{document}